\documentclass[a4paper,10pt]{article}

\RequirePackage{amsthm,amsmath}
\RequirePackage[numbers,square]{natbib}
\RequirePackage[utf8x]{inputenc}
\usepackage{graphicx,graphics,epsfig,latexsym,amssymb,caption,subcaption,float}
\RequirePackage{hyperref}

\usepackage{authblk}

\font\bb=msbm10

\def\R{\hbox{\bb R}}
\def\E{\hbox{\bb E}}
\def\S{\hbox{\bb S}}
\def\N{\hbox{\bb N}}

\def\C{\hbox{\bb C}}

\newtheorem{lemma}{Lemma}
\newtheorem{theorem}{Theorem}
\newtheorem{proposition}{Proposition}
\newtheorem{definition}{Definition}
\newtheorem{corollary}{Corollary}

\newenvironment{Myitemize}%
{\begin{itemize}}%
{\end{itemize}}

\title{Asymptotic Cones of Embedded Singular Spaces}

\author[*]{Xiang Sun}
\author[*,**]{Jean-Marie Morvan}
\affil[*]{King Abdullah University of Science and Technology, KSA}
\affil[**]{Institut Camille Jordan, UMR 5208 CNRS, Universit\'{e} Lyon I, France}

\date{January, 2015}

\begin{document}

\maketitle

\begin{abstract}{
We  use  geometric measure theory to introduce the notion of {\it asymptotic cones} associated with a singular subspace of a Riemannian manifold. This  extends the classical notion of {\it asymptotic directions} usually defined on smooth submanifolds. We get a simple expression of these cones for polyhedra in $\E^3$, as well as convergence and approximation theorems. In particular, if a sequence of singular spaces tends to a smooth submanifold, the corresponding sequence of asymptotic  cones tends to the asymptotic cone of the smooth one for a suitable distance function. Moreover, we apply these results to approximate the asymptotic lines of a smooth surface when the surface is approximated by a triangulation.}
\end{abstract}

\section*{Introduction}

In the past decades, there has been a growing interest in defining geometric invariants describing singular spaces \cite{bernig2003courbures}, \cite{bernig2006curvature}, \cite{cheeger1984curvature}, \cite{cheeger1984kinematic}, \cite{cohen2003restricted}, \cite{cohen2006second},\cite{dai2006accuracy}, \cite{dai2008variational}, \cite{federer1959curvature}, \cite{federer1969geometric}, \cite{fu1993convergence}, \cite{gu2008computational}, \cite{morvan2008generalized}, \cite{zahle1986integral}. Such invariants are generally subject to two assignments:
\begin{enumerate}
\item \label{C1}
They must fit with the classical invariants when the underlying set is a smooth manifold or submanifold.
\item \label{C2}
They must satisfy some continuity conditions. For instance, if a sequence of singular spaces tends (for a suitable topology) to a smooth space, then the invariants defined on the singular spaces also tend to the smooth ones.
\end{enumerate}

After the length, the area and the volume, the most popular smooth geometric invariants are  (sectional, Ricci, scalar, principal) curvatures, which are smooth functions (or tensors) defined on a (smooth) (sub)manifold and its tangent space. A classical approach to extending these curvatures to singular spaces ${\cal W}$ of a Riemannian manifold $M$, is to replace functions by measures on $M$. These measures are defined by integrating invariant differential forms over Borel subsets. As an example,  Lipschitz-Killing curvature measures for singular spaces of a Riemannian manifold can be defined as follows:

\begin{itemize}
\item
The first step is to generalize the unit normal bundle of a smooth submanifold to singular spaces. This has been done for convex subsets and for subsets of positive reach. More generally, the theory of the normal cycle \cite{wintgen1982normal} allows us to define an integral current on a large class of singular subsets, called “geometric subsets” of a Riemannian manifold, generalizing the unit normal bundle.
\item
The second step consists of defining standard differential forms on the tangent bundle of the ambient space.
\item
Finally, integrating these differential forms on the normal cycle,  builds invariant measures, satisfying the two assignments described above.
\end{itemize}

We remark, however, that this building  does not allow us to extract characteristic subsets induced by a couple $({\cal W},M)$  (where ${\cal W}$ is a singular space of a Riemannian manifold $M$) as asymptotic directions and principal directions.

In the framework of geometric measure theory, more specifically using the normal cycle theory,  we use the formalism introduced in \cite{cohen2003restricted},\cite{cohen2006second} and \cite{morvan2008generalized} to propose an extension of the definition of the classical {\it asymptotic directions} to a large class of singular spaces ${\cal W}$. Instead of building a new ``curvature measure" on any couple $({\cal W},M)$, we associate a map
 that assigns  a {\it cone} (called an {\it asymptotic cone}) of vector fields on $M$ to each Borel subset of $M$. In particular, if $M=\E^N$, we can reduce the target of this map  to the set of quadratic cones of $\E^N$. Moreover, choosing a fixed scalar $r>0$, we associate a fie	ld of cones leaving in the tangent bundle $T\E^N$ with such a couple $({\cal W},\E^N)$.
The two assignments described above are therefore satisfied.
\begin{itemize}
\item
The cones associated with a Borel subset reduced to a single point of a smooth surface lying in $\E^3$ are nothing but the product of the classical asymptotic direction(s) by the normal line of the surface at this point.
\item
If ${\cal W}_k$ is a sequence of singular subsets (admitting asymptotic cones), which tends to a smooth submanifold $M$, then the sequence of asymptotic cones associated with ${\cal W}_k$ built over a  ``regular" Borel subset $B$ of $M$ tends to the cone associated with $M$ built over the same $B$.

\end{itemize}

We give some applications of this building up for surfaces in $\E^3$.
\begin{itemize}
\item
If a plane field is defined on a singular surface of $\E^3$, we can define the asymptotic directions at a point $m$ of the surface associated with a Borel set (for instance a ball centered at $m$) as the intersection of the asymptotic cone defined at $m$,  with the plane defined at $m$, to get two directions than can be called {\it asymptotic directions}. By integrating these directions, we get {\it asymptotic lines}
depending obviously on the plane field.
\item
Using a triangulated approximation  of a smooth surface $W^2$ in $\E^3$, we build an approximation of its asymptotic directions at some fixed points by considering the intersection of the asymptotic cones of the triangulation with the planes spanned by the faces of the triangles. We can deduce  an approximation of the  asymptotic lines of $W^2$ if the triangulation is a ``good" approximation of the smooth surface.
\end{itemize}

%

Finally, we mention that this approach can be extended to other characteristic subsets or subspaces, which could be the subject of  future researches.

The article is organized as follows:

After a summary of the main notations pertaining to Riemannian submanifolds, the first section  begins  with smooth surfaces in $\E^3$.

Since we  define a new geometric invariant in quite a general context, our framework is an (oriented ${\cal C}^{\infty}$) Riemannian manifold $M$ (of finite dimension). After a review of the geometry of the tangent bundle of  $M$ (see \cite{grifone1972structureI} for instance), the second section introduces the so-called {\it asymptotic form} defined in the tangent bundle of $M$, easily derived from the fundamental form \cite{cohen2006second}, \cite{morvan2008generalized}.

The third section uses the theory of the normal cycle \cite{fu1988monge}, \cite{wintgen1982normal} to introduce the notion of {\it asymptotic measures} defined in the tangent bundle of $M$, associated with a large class of singular subsets $\cal G$  (called {\it geometric subsets} in \cite{fu1988monge}).  These measures are derived from the asymptotic form defined in Section \ref{ASFO}.
We describe them explicitly in  classical situations (when $\cal G$ is smooth or a polyhedron).
Then, we construct the  {\it asymptotic cones}. In particular, we show how this construction generalizes the asymptotic directions defined on a smooth surface in $\E^3$, fulfilling    assignment \ref{C1}  above. In particular,
we give an explicit expression of  the asymptotic cones associated with a $2$-polyhedron in $\E^3$.

The fourth section  deals with a general theorem of convergence of asymptotic cones. This is the justification of our definition, and fulfills assignment \ref{C2}. In particular,  we give explicit assumptions in terms of the fatness of the sequence of polyhedra, implying the convergence of the sequence of cones for a suitable pseudo-metric.
The last section presents two applications. First, we give a method to build asymptotic lines on a triangulation. Then, we give a method to approximate the asymptotic lines of a smooth surface approximated by a polyhedron. We test our method on various smooth or discrete surfaces.

{\bf Some notation -}
For details on the theory of smooth Riemannian submanifolds, the reader may consult for instance, \cite{chen1973geometry} or \cite{do1992riemannian}. We use the following notation.
Let $W^n$ be an $n$-dimensional closed (oriented) smooth submanifold embedded in a smooth $N$-dimensional (oriented) Riemannian manifold $(M^N,<.,.>)$.
The manifold  $W^n$ inherits a Riemannian structure by pulling back  $<.,.>$ by the embedding, which we still denote by $<.,.>$.  We denote by
 $TW^n\overset{\pi_{W^n}}{\longrightarrow}W^n$
  ({\it resp.} $TM^N\overset{\pi_{M^N}}{\longrightarrow}M^N$)
  the tangent bundle of $W^n$ ({\it resp.} $M^N$) and
   $\Xi (W^n)$ ({\it resp.} $\Xi(M^N)$) the space of tangent vector fields over $W^n$ ({\it resp.} $M^N$).  We denote by $\tilde{\nabla}$ ({\it resp.} $\nabla$)  the Levi-Civita connection on $(M^N,<.,.>)$ ({\it resp.} $W^n$). We denote by $T^{\perp}W^n\overset{\pi_{W^n}}{\longrightarrow}W^n$ the normal bundle of the submanifold $W^n$. The {\it second fundamental form} of the submanifold $W^n$ is the symmetric vector-valued $(2,0)$-tensor
$$h~: TW^n \times TW^n \to T^{\perp}W^n$$
defined as follows:
$$\forall x \in TW^n, \forall y \in TW^n, h(x,y)= \tilde{\nabla}_xy- \nabla_xy.$$
Let $m$ be a point of $W^n$. The isotropic cone
\begin{equation}\label{IC}
{\bf C}^{W^n}_m=\{x\in T_mW^n~: h_m(x,x)=0\}
\end{equation}
 of $h_m$ is classically called the {\it asymptotic cone} of $W^n$ at $m$.
 For any $\xi_m \in T_m^{\perp}W^n$, the eigenvalues of $<h_m(.,.),\xi_m>$ are  the {\it principal curvatures} of $W^n$ at $m$ in the direction $\xi_m$.

\section{The case of smooth surfaces in $\E^3$} \label{SURFACES}
We first restrict our attention to smooth closed (oriented) surfaces $W^2$ embedded in the (oriented) Euclidean space $(\E^3,<.,.>)$ bounding a domain $D$. Let $\xi$ be the normal vector field compatible with these orientations. The second fundamental form of $W^2$ can now be identified with the tensor
$<h(.,.),\xi>$ taking its values in ${\cal C}^{\infty}(W^2)$. We denote by  $\lambda_{1_m},\lambda_{2_m}$ the principal curvatures of $W^2$ at the point $m$, that is, the eigenvalues of $h_m$, by $G$ the Gauss curvature of $W^2$, that is, the determinant of $h$ and by $H$ its mean curvature, that is, its trace (in an orthonormal frame). In a  frame of principal vectors $(e_{1_m}, e_{2_m})$  at $m$ (that is, eigenvectors of $h_m$) the matrix of $h_m$ is
$$\begin{pmatrix}
\lambda_{1_m} & 0 \\
0 & \lambda_{2_m}
\end{pmatrix}.
$$
At each point $m \in W^2$ with negative Gauss curvature, the asymptotic cone  ${\bf C}_{m}^{W^2}$ is the union of two lines. Integrating the corresponding vector fields gives rise to  foliations of $W^2$ by the so-called {\it asymptotic curves}. By definition, at each point, the principal normal directions of these curves (considered as curves in $\E^3$)  are tangent to $W^2$.

 Using measure theory, the goal of this paper is to define and study analogous cones  associated with a   (regular or) singular subspace of a Riemannian manifold $M^N$, lying above any Borel subset of $M^N$.

Let us begin by explaining how we define such cones over any Borel subset of $\E^3$ in the regular case; that is, when the subset is a (compact) domain bounded by a smooth surface $W^2$. Let $T_{W^2}\E^3$ be the tangent bundle of $\E^3$ restricted to $W^2$. If $x$ is any vector field on $T_{W^2}\E^3$,
we build a  signed  measure $\Phi_{W^2}^x$ as follows: For any Borel subset $B$ of $\E^3$, we write
\begin{equation}\label{SURFACE2}
\Phi_{D}^x(B) = \int_{B\cap W^2}h_m(pr_{T_mW^2}x,pr_{T_mW^2}x)dm,
\end{equation}
where $pr_{T_mW^2}$ denotes the orthogonal projection on $T_mW^2$ and $dm$ is the Lebesgue measure on $\E^3$.
Let us now fix $B$ and consider the map
\begin{equation}
x \mapsto\Phi_{D}^x(B),
\end{equation}
where $x$  runs over the (huge) space of
vector fields  $\Xi(\E^3)|_{W^2}$. This map is quadratic in $x$. If we  force $x$ to be a constant vector field, then we  get a quadratic form (that we still denote by $\Phi^{\bullet}_{D}(B)$) on $\E^{3}$.
This quadratic form has generically three eigenvalues, $\lambda_1(B),\lambda_2(B),\lambda_3(B)$, that we call the {\it principal curvatures} of $B$. The corresponding eigenvectors are called the {\it principal vectors} of $B$, and the matrix of
 $\Phi^{\bullet}_{W^2}(B)$ in this frame is
 $$\begin{pmatrix}
 \lambda_1(B) & 0 & 0 \\
 0 & \lambda_1(B) & 0 \\
 0 & 0 & \lambda_3(B)
 \end{pmatrix}.
 $$
From this construction, we also introduce the {\it isotropic cone} associated with $\Phi_{D}^{\bullet}(B)$:
\begin{equation}
{\cal C}^{{\rm par}, D}_B = \{x \in \E^3; \Phi_{D}^x(B)=0\},
\end{equation}
(the notation coming from the fact that we restrict our isotropic cone to constant - that is, parallel - vector fields in $\E^3$).
We call it the {\it asymptotic cone} of $B$ (with respect to $D$).


To clarify  that this construction is linked with the classical pointwise situation, suppose  that $B$ is reduced to a point $\{m\}\in W^2$.
 If $y$ is a constant vector field such that $y_m\in T_mW^2$, then
\begin{equation}
\Phi_{D}^y(\{m\})=h_m(y_m,y_m).
\end{equation}

If $z$ is a constant vector field such that $z_m=\xi_m$, then
\begin{equation}
\Phi_{D}^z(\{m\})=0.
\end{equation}

This implies that in the frame $(e_{1_m},e_{2_m},  \xi_m)$, the matrix of $\Phi_{D}^{\bullet}(\{m\})$ is
\begin{equation}
\begin{pmatrix}
\lambda_{1_m}  & 0 & 0 \\
0 &\lambda_{2_m} &0\\
0 & 0 & 0
\end{pmatrix}.
\end{equation}
Consequently, the asymptotic cone ${\cal C}^{{\rm par}, D}_m$ is nothing but the cone spanned by the normal $\xi_m$ and ${\bf C}^{W^2}_m$.


This construction has some advantages:
One can define the asymptotic cones  at different scales by scaling the Borel sets (for instance, by taking balls of radius $\frac{1}{k}$ as Borel subsets).
Generically, we get three geometric invariants, $\lambda_1(B),\lambda_2(B),\lambda_3(B)$, instead of two, and a  two-dimensional cone instead of the union of two lines.
Moreover, another important  advantage of replacing functions by measures, is that this framework can be used for a large class of singular spaces (for instance, polyhedra,  algebraic subsets, subanalytic subsets) of any codimension in any Riemannian manifold, as long as we can extend the notion of normal space.
For instance, if one replaces the smooth surface by a polyhedron $P$ bounding a domain ${\cal D}$, we will get the following explicit simple expression approximating the cone ${\cal C}_{B}^{\text{par},{\cal D}}$ (see \ref{JQQZ56}, and also \ref{POLY23}  for the exact formulas)~:

$$
 \Phi_{\cal D}^x(B) \sim \sum_{e \in {\bf E}}l(e \cap B)
 \angle(e)<x,e^{-}>^{2},
$$
and
$$
{\cal C}_{B}^{\text{par},{\cal D}} \sim \{x \in \E^{3}: \sum_{e \in {\bf E}}l(e \cap B)
 \angle(e)<x,e^{-}>^{2}=0\},
$$
(see Section \ref{NMA} for the notation).
This is why the theory of the normal cycle, extensively studied over the last decades \cite{fu1988monge}, \cite{wintgen1982normal}, will be our framework.
In the following we describe the construction of asymptotic cones in such a large context.

\section{Asymptotic forms} \label{ASFO}
To be self contained, we begin with a summary of the geometry of the tangent bundle of an (oriented) $N$-dimensional smooth Riemannian manifold  $(M^N, <.,.>)$.
The reader may consult \cite{grifone1972structureI}, \cite{grifone1972structureII} and \cite{cohen2006second} for details.
We denote by $TTM^N\overset{\pi_{TM^N}}{\rightarrow}TM^N$ the tangent bundle of the manifold $TM^N$. As usual, we consider the exact sequence of vector bundles:
\begin{equation}\label{seq::exact}
0\longrightarrow TM^N\times_{M^N} TM^N \overset{i}{\longrightarrow} TTM^N \overset{j}{\longrightarrow} TM^N \times_{M^N} TM^N \longrightarrow 0,
\end{equation}
where $i$ denotes the natural injection defined by
\begin{equation}
i(u_1,u_2) = \frac{d}{dt}(u_1+tu_2)\Big| _{t = 0}
\end{equation}
and
\begin{equation}
j = (\pi_{TM^N},d\pi_{M^N}).
\end{equation}
The {\it vertical bundle} of $M$ is the subbundle $V(M^N)=\ker j$ of $TTM^N$. The morphism $i$ induces an isomorphism:
$$\overline{i} : TM^N\times_{M^N} TM^N  \to V(M^N).$$
If $m \in M^N$ and $x\in T_mM^N$, the {\it vertical lift} of $z\in T_mM^N$ at $x$ is the vector $z^v = i(x,z)$.
The morphism $J = i\circ j$ is an almost tangent structure on $M^N$  ($J^2 = 0$) and $V(M^N) = \ker J$. Let
\begin{equation}
\delta :TM^N \longrightarrow TM^N \times_{M^N} TM^N
\end{equation}
 be the canonical vector field defined by
$\delta (x)=(x,x)$
and let $C: TM^N\longrightarrow V(M^N)$ be the vertical vector field associated with the (global) one-parameter group of homotheties with positive ratio, acting on the fibers of $TM^N$. We have
$C=i\circ\delta$.
%
 We write
$\eta = pr_2 \circ \overline{i}^{-1}$,
where  $pr_2$ denotes the projection on the second factor of
$TM^N\times_{M^N} TM^N$.
Since $M^N$ is endowed with a Riemannian metric $<\cdot,\cdot>$ and its Levi-Civita connection, we can build the corresponding right splitting
\begin{equation}
\gamma: TM^N \times_{M^N} TM^N \longrightarrow TTM^N
\end{equation}
of the exact sequence \ref{seq::exact},
(satisfying $j\circ\gamma = {\rm Id}_{TM^N\times_{M^N} TM^N}$).
Let $m \in M^N$ and $x\in T_mM^N$. {\it The horizontal lift} of $z\in T_mM^N$ at $x$ is the vector $z^h = \gamma(x,z)$.
We denote $H_x(M^N) = {\rm Im}(\gamma(x,\cdot))$, from which we construct the horizontal bundle $H(M^N)$ such that, for all $x\in T_mM^N$
\begin{equation}
T_x TM^N = V_x (M^N)\oplus H_x (M^N).
\end{equation}
 We denote by ${\cal V}(M^N)$ ({\it resp.} ${\cal H}(M^N)$) the space of vertical ({\it resp.} horizontal) vector fields.
We denote by $\mathtt{h}: TTM^N\longrightarrow H(M^N)$ the horizontal projection, and by $\mathtt{v}: TTM^N\longrightarrow V(M^N)$ the vertical projection. We remark that $\mathtt{h} = \gamma\circ j$. 
The morphism
$$K=\eta \circ \mathtt{v} :TTM^N \to TM^N$$
is the {\it connector} associated with the Levi-Civita connection. At every point
$x \in T_mM^N$, the morphism
$$(d\pi_M \times K)_x : T_xTM \to T_mM^N \times T_mM^N$$
is an isomorphism that identifies $V_x(M^N)$ with $T_mM^N$  and $H_x(M^N)$ with $T_mM^N$.
The bundle $TTM^N\overset{\pi_{TM^N}}{\longrightarrow} TM^N$ is canonically endowed with the {\it Sasaki metric} $<\cdot,\cdot>$ defined by the following conditions:
\begin{equation}
\left\{
\begin{aligned}
&V(M^N) \text{ and } H(M^N) \text{ are orthogonal,}\\
&i \text{ is an isometry,}\\
&\gamma \text{ is an isometry.}
\end{aligned}\right.
\end{equation}
If
$$\alpha : M^N \to M^N \times M^N$$
is the diagonal map defined by $\alpha(m)=(m,m)$, then
for every $x \in TM^N$,
$$(d\pi_M \times K)^{-1} \circ d\alpha (x)= x^v \oplus x^h.$$
Finally, the bundle $TTM^N\overset{\pi_{TM^N}}{\longrightarrow} TM^N$ is also endowed with an almost complex structure $F$ ($F^2 = -\text{Id}$), defined by the following conditions
\begin{equation}
\left\{
\begin{aligned}
&FJ=\mathtt{h}\\
&F \mathtt{h}=-J.
\end{aligned}
\right.
\end{equation}
Therefore $F|_{V}:V(M^N)\longrightarrow H(M^N)$ and $F|_{H}:H(M^N)\longrightarrow V(M^N)$ are isometries.
In this Riemannian context, we give the following definition:

\begin{definition}
\begin{enumerate}
\item
 The vector valued $(N-1)$-form on $TM^N$  defined for each $X\in {\cal H}(M^N)$ by
\begin{equation}
\mathbf{h}^X=[*_{\rm Hodge}(FC\wedge X)]\wedge FX
\end{equation}
is called the {\rm asymptotic  $(N-1)$-form} on $TM^N$.
\item
The vector valued $(N-1)$-form on $M^N$  defined for each $x \in \Xi(M^N)$  by
\begin{equation}
\mathbf{h}^{x}= [*_{\rm Hodge}(FC\wedge x^h)]\wedge x^v
\end{equation}
is called the {\rm asymptotic  $(N-1)$-form} on $M^N$.
\end{enumerate}
\end{definition}

In this definition, $*_{\rm Hodge}$ denotes the Hodge duality on $H_x(M^N)$ for each $x\in TM^N$. (The introduction of the Hodge operator in the definition of the generalized second fundamental form can be found in \cite{cohen2008stability} when the ambient space is Euclidean. We adapt it here to the general Riemannian situation. It is equivalent to the initial definition given in  \cite{cohen2006second}). Using the identification of vector fields and $1$-forms induced by the Riemannian structure, $*_{\rm Hodge}(FC\wedge X)$ is a $(N-2)$-form on ${\cal H}(M^N)$. On the other hand, $FX$ is (identified with) a $1$-form, null on $H(M^N)$ and acting on $V(M^N)$, and  $*_{\rm Hodge}(FC\wedge X)$ is (identified with) an $(N-2)$-form null on $V(M^N)$ and acting on  $H(M^N)$.

\section{Normal cycles, asymptotic  measures, asymptotic cones} \label{NMA}
\subsection{Currents and normal cycles of singular spaces}
Let ${\cal D}_l(TM^N)$ be the space of {\it $l$-currents} of $TM^N$ ($0 \leq l \leq 2N$); that is,  the topological dual of  the space  ${\cal D}^l(TM^N)$ of {\it $l$-differential forms} with compact support on $TM^N$, endowed with the topology of uniform convergence on any compact subset, of all partial derivatives of any order.  The duality bracket will be still denoted by $<\cdot,\cdot>$ if no confusion is possible.
The space ${\cal D}_l(TM^N)$ is naturally endowed with the weak topology: if $(C_k)_{k \in \N}$ is  a sequence of $l$-currents of $TM^N$ and if $C$ is a $l$-current of $TM^N$, then
\begin{equation}
\lim_{k\rightarrow\infty} C_k = C \Longleftrightarrow \forall \omega \in {\cal D}^l(TM^N), \lim_{k\rightarrow\infty}< C_k,\omega >=< C,\omega >.
\end{equation}
An $l$-current is {\it rectifiable} if it is associated with a rectifiable subset (see \cite{morvan2008generalized} for details). An $l$-current is {\it integral} if it is rectifiable and its boundary is rectifiable.

%
%


When it exists, the {\it normal cycle} of a (compact singular) subset ${\cal W}$ of a Riemannian manifold $M^N$ is a closed integral current
 ${\bf N}({\cal W}) \in {\cal D}_{N-1}(TM^N)$, which is Legendrian for the symplectic structure on $TM^N$ dual to the canonical one on $T^{*}M^N$ in the duality defined by the metric. The normal cycle is the direct generalization of the unit normal bundle of a smooth submanifold.  Its formal definition was given in \cite{fu1988monge}. Although the normal cycle cannot be defined on any compact  subset of $M^{N}$, it exists for a large class of  subsets, as convex subsets,  polyhedra, subsets of positive reach, subanalytic subsets for instance. Following \cite{fu1988monge}, any compact subset ${\cal G}$ of $M^N$ such that ${\bf N}({\cal G})$ exists is
said to be {\it geometric}, and ${\bf N}({\cal G})$ is called its normal cycle. One of the  main properties of the normal cycle for our purpose is its additivity
\cite{fu1988monge}:
\begin{proposition}\label{proposition::ADD}
If ${\cal G}_1$ and ${\cal G}_2$   are geometric, then
${\cal G}_1 \cup {\cal G}_2$ and ${\cal G}_1 \cap {\cal G}_2$ are geometric
and
\begin{equation}\label{equation::ADD}
{\bf N}({\cal G}_1  \cup {\cal G}_2)= {\bf N}({\cal G}_1) + {\bf N}({\cal G}_2) - {\bf N}({\cal G}_1 \cap {\cal G}_2).
\end{equation}
\end{proposition}

Here are some classical examples:
\begin{enumerate}
\item \label{E1}
The normal cycle of a smooth submanifold of a Riemannian manifold is the closed current associated with its unit normal bundle.
\item \label{E2}
If $D$ is a compact domain whose boundary is a smooth hypersurface, then its normal cycle is the closed current associated with its outward unit normal vector field.
\item \label{E3}
If $C$ is a convex body, then its normal cycle is the closed current associated with the oriented set
$$\{(m, \xi) :
 m \in \partial C, \xi \in \E^3 , ||\xi|| = 1, \forall z \in C, <\xi,\overrightarrow{mz}> \leq 0\}.$$
\item \label{E4}
The normal cycle of a polyhedron of $\E^N$ can be computed  by applying \ref{equation::ADD} to a decomposition of the polyhedron into (convex) simplices and using \ref{E3}.
\end{enumerate}

\subsection{Asymptotic measures, asymptotic cones} \label{ASMASC}
Let us  now define an  {\it asymptotic (signed) Radon  measure} on $M^N$ ({\it resp.} $TM^N$) associated with a geometric subset. We denote by  ${\cal B}_{M^N}$ ({\it resp.} ${\cal B}_{TM^N}$)  the class of Borel subsets of $M^N$  ({\it resp.} $TM^N$) with compact closure.

\begin{definition}
Let ${\cal G}$ be a geometric subset of $M^N$.
\begin{Myitemize}
\item
The {\rm asymptotic  measure} defined on $TM^N$, associated with ${\cal G}$ and  $X\in{\cal H}(M)$ is the map
\begin{equation}
{\Phi}_{\cal G}^X: {\cal B}_{TM^N}\longrightarrow \R
\end{equation}
defined as follows:
\begin{equation}
\forall B \in {{\cal B}_{TM^N}}, {\Phi}_{\cal G}^X({B}) = <{\bf N}({\cal G}),\chi_{{B}}\mathbf{h}^X>.
\end{equation}

\item
The {\rm asymptotic  measure} defined on $M^N$, associated with ${\cal G}$ and $x\in \Xi(M^N)$ is the map
 ${\Phi}_{\cal G}^x$ defined as follows:
\begin{equation}\label{GENEX}
\forall B\in {\cal B}_{M^N}, {\Phi}_{\cal G}^{x}(B) = <{\bf N}({\cal G}),\chi_{\pi_{M}^{-1}{(B)}}\mathbf{h}^{x}>.
\end{equation}
\end{Myitemize}
\end{definition}

If ${\cal G}$ and $B$ are fixed, the map
$$x \to <{\bf N}({\cal G}),\chi_{\pi_{M}^{-1}{(B)}}\mathbf{h}^{x}>$$
is quadratic, inducing its isotropic cone. This remark leads to the following definition:

\begin{definition}
Let ${\cal G}$ be a geometric subset of $M^N$. With any Borel subset $B  \in {\cal B}_{M^N}$, we associate the cone
\begin{equation}
{\cal C}_{B}^{\cal G}=\{x \in \Xi (M^N):\Phi_{\cal G}^{x}(B)=0\}
\end{equation}
and the cone
\begin{equation}
{\cal C}_{B}^{{\rm par},{\cal G}}=\{x \in \Xi (M^N) :  x \mbox{ \rm parallel}, \Phi_{\cal G}^{x}(B)=0\},
\end{equation}
which we call the {\rm asymptotic cones} associated with $B$.
\end{definition}

In many applications, and for simplicity, it is easier to  consider the cone of {\it parallel} vector fields, identifying a parallel vector field with its value at any point $m \in B$.
Obviously, ${\cal C}_{B}^{\text{par},{\cal G}} \subset   {\cal C}_{B}^{\cal G}$. We remark, however, that this new definition can be quite restrictive depending on the geometry of $M^N$. For instance, if $M^N$ has non-zero constant sectional curvature, the only parallel vector field is the null vector field. In contrast, if $M^N=\E^N$, the space of parallel vector fields is the space of constant vector fields, which can be identified with $\E^N$, which is much easier to manipulate. We deduce easily from \ref{GENEX} explicit expressions of these curvature measures in some particular cases.

\subsubsection{The case of smooth submanifolds}
Let $W^n$ be a (compact smooth) submanifold (with or without a boundary) embedded in $M^N$ and  $X\in {\cal H}(TM^N)$. Since the normal cycle of a smooth submanifold is its unit normal bundle,  we deduce from \cite{cohen2006second} or \cite{morvan2008generalized} (Corollary $16$ page $215$) that for any $B \in {\cal B}_{TM^N}$,
\begin{equation}
{\Phi}_{W^{n}}^X(B) = \int_{ST^\perp W^n\cap B}h^\xi(\text{pr}_{TW^n}d\pi_M(X),\text{pr}_{TW^n}d\pi_M(X))d\xi dv,
\end{equation}
where $h^{\xi}$ denotes the second fundamental form of $W^n$ in the direction of the unit vector $\xi$, $ST^\perp W^n$ denotes the unit normal bundle of $W^n$ and $\text{pr}_{TW^n}$ denotes the orthogonal projection onto the tangent bundle $TW^n$.
In particular, let $W^{N-1}$ be a (smooth oriented) hypersurface of $M^N$  bounding  a domain $D$.
(This assumption is not restrictive in our case, since our results are local. It allows to simplify some technical points by considering ``only one side"  of the normal cycles (the one corresponding to the outward unit normals). we have, for any $B$ in ${\cal B}_{M^N}$,
\begin{equation}
 {\Phi}_D^X(B) = \int_{W^{N-1}\cap B}h(\text{pr}_{TW^{N-1}}d\pi_M(X),\text{pr}_{TW^{N-1}}d\pi_M(X)) dv,
\end{equation}
where $\xi$ is the outward (with respect to $D$) unit normal vector field of $W^{N-1}$, and $h$ is the second fundamental form of $W^{N-1}$ in the direction $\xi$. We have then a correct generalization of \ref{SURFACE2}. Consequently,
\begin{itemize}
\item
If $B$ is reduced to a point $m$,
\begin{equation}\label{LJRE}
{\cal C}_m^{W^n}= \{ x\in T_m M^N: h_m(pr_{TW^n}x, pr_{TW^n}x)=0\}.
\end{equation}
We deduce that  ${\cal C}_{m}^{W^n}$ is the cone spanned by ${\bf C}_{m}^{W^n}$ and $T_m^\perp W^n$; that is, we have the direct generalization of the corresponding cone defined for surfaces in $\E^3$ in Section \ref{SURFACES}.

\item
 If $W^{N-1}$ is a (smooth-oriented) hypersurface of $\E^N$ bounding a domain $D$,
then
\begin{equation}\label{278}
{\cal C}_{B}^D=\{x \in TM^N : \int_{W^{N-1} \cap B}h(pr_{TW^{N-1}}x,pr_{TW^{N-1}}x)dv=0\},
\end{equation}
and
\begin{equation}\label{279}
{\cal C}_{B}^{\text{par},D}=\{x \in \E^{N}:\int_{W^{N-1}\cap B}h(pr_{TW^{N-1}}x,pr_{TW^{N-1}}x)dv=0\}.
\end{equation}
\end{itemize}
\subsubsection{The case of  polyhedra}
We will extend \ref{278} and \ref{279} to polyhedra. Let ${\cal D}$ be a domain  in $\E^N$ bounded by a   $(N-1)$-dimensional  polyhedron $P^{N-1}$. For $X \in {\cal H}(P^{N})$,
we can  evaluate $\Phi_{\cal D}^X$ above each simplex. In particular, if $\sigma^{N-2}$ is a   $(N-2)$-simplex, the support of $N({\cal D})|_{\sigma^{N-2}}$ is the product $\sigma^{N-2} \times C_\sigma$, where $C$ is a portion of circle. Let $(e_1,..., e_{N-2})$ be an orthonormal frame field tangent to $\sigma^{N-2}$. Any point of $\sigma^{N-2} \times C_\sigma$ is a couple $(m,e_{N-1})$, where $m$ is a point of $\sigma^{N-2}$ and $e_{N-1}$ is a unit vector orthogonal to $\sigma^{N-2}$. With these notations, we deduce from \cite{cohen2006second} or \cite{morvan2008generalized} (Theorem $72$ page $216$) that  for
 any $B  \in {\cal B}_{M^N}$,
\begin{equation}\label{POLY1}
\Phi_{\cal D}^X(B) =\sum_{\sigma^{N-2}\subset\partial P^N}\int_{(\sigma^{N-2}\cap B)\times C}< X, e^h_{(N-1)} >^2.
\end{equation}
We also deduce that for any $B \in {\cal B}_{\E^{N}}$,
\begin{equation}\label{POLY2}
{\cal C}_{B}^{\text{par},{\cal D}}=\{x \in \E^{N}: \int_{(\sigma^{N-2}\cap B)\times C}< x^h, e^h_{(N-1)} >^2=0\}.
\end{equation}
In particular, if ${\cal D}$ is a domain of
 $\E^{3}$  bounded by a polyhedron $P^2$, \ref{POLY2} can be reduced to an explicit simple expression: First of all, we identify the (vector) plane  $e^{\perp}$ orthogonal to any (oriented) edge $e$ of $P$  with $\C$, as follows~:  Let $n_1 \in e^{\perp}$ ({\it resp.} $n_2 \in e^{\perp}$) be the unit (oriented) normal to  the faces $f_1$, ({\it resp.} $f_2$) incident to $e$.  Let $e^{+}$ ({\it resp.} $e^{-}$)  be the (oriented) normalized vectors spanning the bisectors of $n_1$ and $n_2$ (so that $(e^{+},e^{-},e)$ is a direct frame of $\E^3$). Any vector $ae^{+}+be^{-}$ of $e^{\perp}$ is now identified with the complex number $a+ib$. An  explicit integration over each term of type $(\sigma^1 \cap B) \times C$ in \ref{POLY2} gives the following expression of the asymptotic measure and asymptotic cone:

 \begin{proposition}
 \begin{enumerate}
 \item
For any $B  \in {\cal B}_{\E^{3}}$ and any {\it constant} vector field $x$  of $\E^{3}$,
\begin{equation} \label{JQQZ5}
 \Phi_{\cal D}^x(B) =\sum_{e \in {\bf E}}\frac{l(e \cap B)}{2}
 \big[(\angle(e)||pr_{e^{\perp}}x||^{2} + \sin (\angle(e)){\cal R}\big((pr_{e^{\perp}}x)^{2}\big)\big],
\end{equation}
where $\bf E$ denotes the set of edges of $P^2$, $\angle(e)$ the angle of the normal to the faces incident to $e$ (being positive if and only if $e$ is convex), $pr_{e^{\perp}}$ the orthogonal projection on $e^{\perp}$, and ${\cal R}(pr_{e^{\perp}}x)^{2}$  the real part of the complex number $(pr_{e^{\perp}}x)^{2}$.
\item
In particular, {\small
\begin{equation}\label{POLY23}
{\cal C}_{B}^{{\rm par},{\cal D}}=\{x \in \E^{3}: \sum_{e \in {\bf E}}\frac{l(e \cap B)}{2}
 \big[(\angle(e)||pr_{e^{\perp}}x||^{2} + \sin (\angle(e)){\cal R}\big((pr_{e^{\perp}}x)^{2}\big]=0\}.
 \end{equation}
 }
 \end{enumerate}
\end{proposition}

We remark that if the $\angle(e)$'s are ``small enough" (this can happen for instance when $P$ approximates ``smoothly" a smooth surface), then $\sin (\angle(e))$ is ``close to" $\angle(e)$ and
\begin{equation} \label{JQQZ55}
 \Phi_{\cal D}^x(B) \sim \sum_{e \in {\bf E}}l(e \cap B)\angle(e)<x,e^{-}>^{2}.
\end{equation}
And then,
\begin{equation} \label{JQQZ56}
{\cal C}_{B}^{{\rm par},{\cal D}} \sim \{x \in \E^{3}: \sum_{e \in {\bf E}}l(e \cap B)
 \angle(e)<x,e^{-}>^{2}=0\}.
\end{equation}

%
%

After  choosing  a scale $r$, we construct from the previous construction, a {\it cone subbundle} of $TM^N$ associated with a geometric subset of $M^N$.
Let us denote by $B(m,r)$ the ball of radius $r$, centered at $m\in M^N$. With each point $m\in M^N$, and for a fixed (small enough) real number $r>0$, we associate the cone
${\cal C}_{B(m,r)}^{\text{par},{\cal G}}$.

\begin{definition}
We call $\cup_{m\in\E^N}{\cal C}_{B(m,r)}^{{\rm par},{\cal G}}$ the {\rm cone subbundle} of $T\E^N$ at scale $r$ associated with ${\cal G}$.
\end{definition}

\begin{figure}[H]
	\centering
	\begin{subfigure}[t]{0.31\textwidth}
                \includegraphics[width=\textwidth]{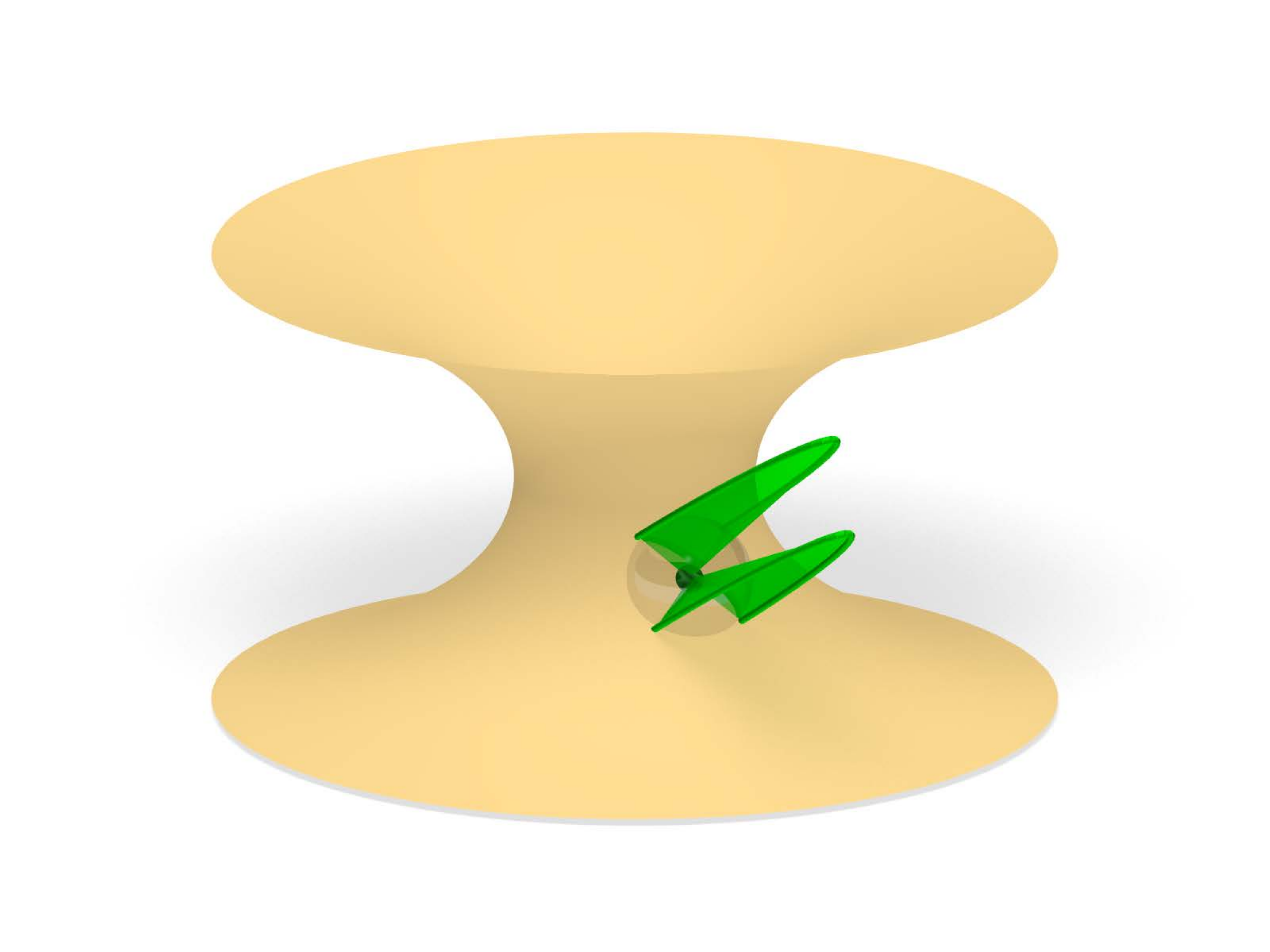}
      	 \subcaption{}\label{subfconeleft}
      \end{subfigure}
	\begin{subfigure}[t]{0.31\textwidth}
        \includegraphics[width=\textwidth]{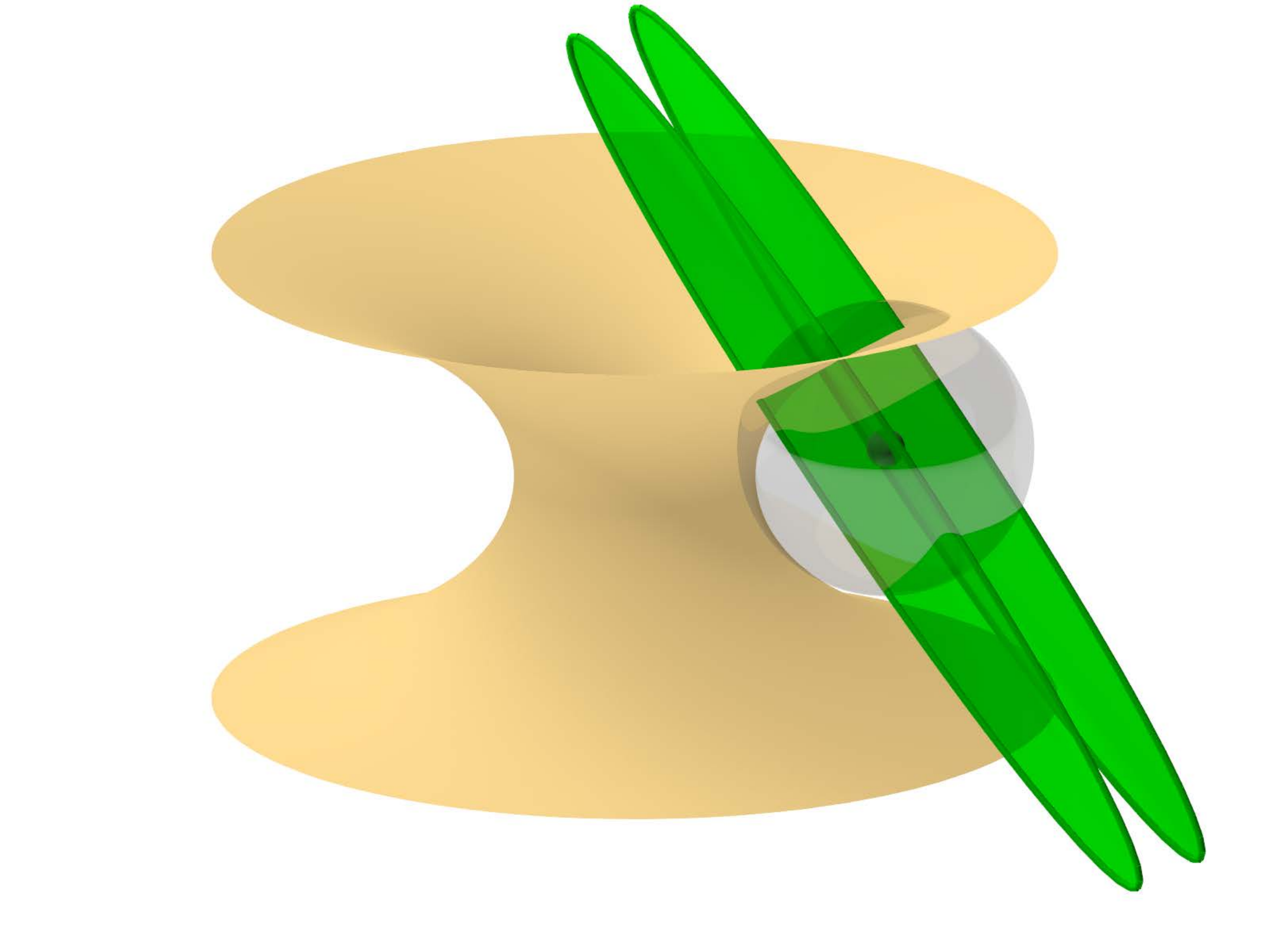}
        \subcaption{}\label{subfconemiddle}
    \end{subfigure}
    \begin{subfigure}[t]{0.31\textwidth}
        \includegraphics[width=\textwidth]{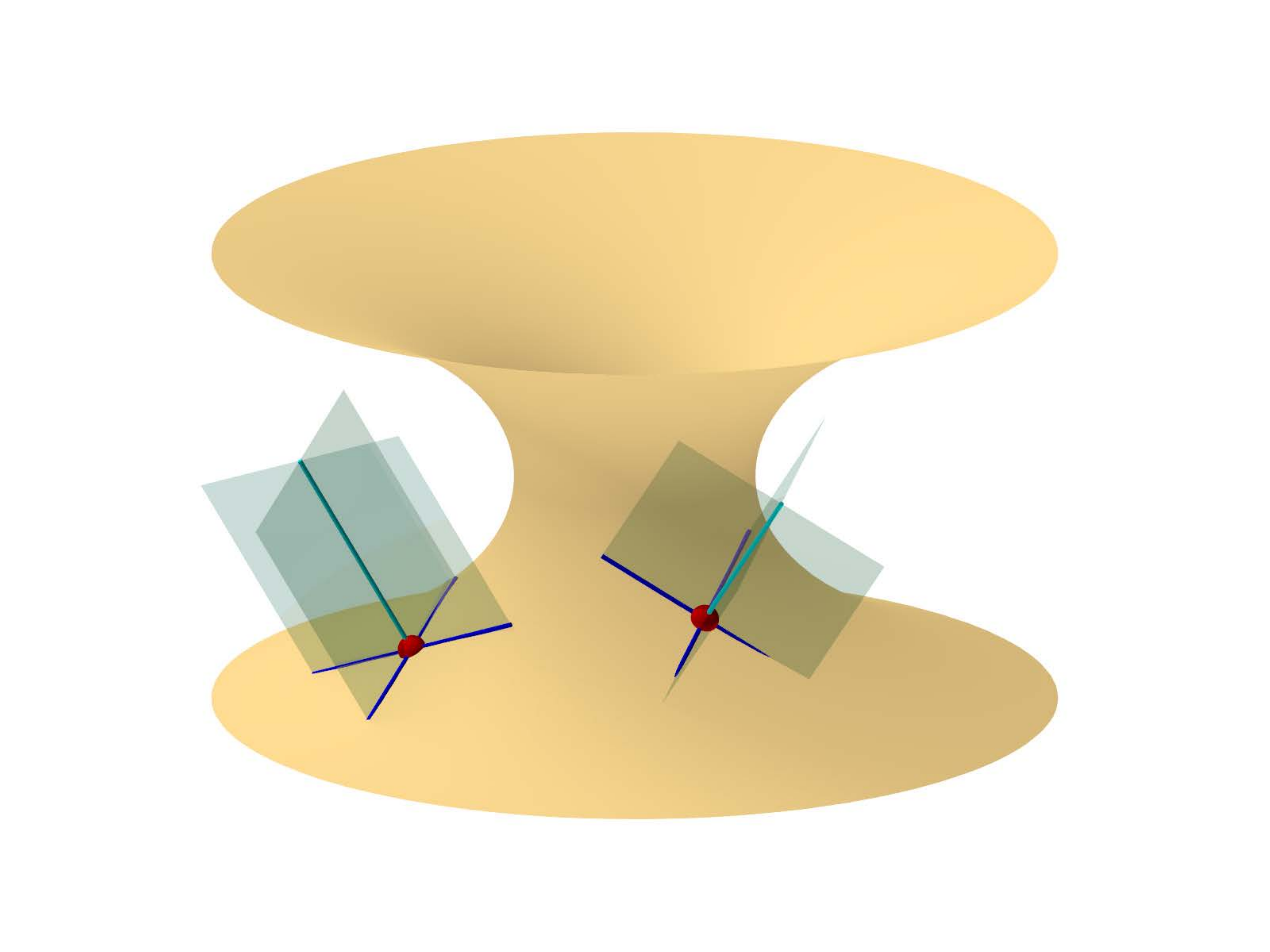}
        \subcaption{}\label{subfconeright}
    \end{subfigure}
    \caption{Some asymptotic cones built on a smooth surface (here on a portion of a catenoid) (a) An asymptotic  cone in green whose vertex is at the center of its corresponding Borel set (a transparent ball) (b) The center of the ball may be out of the surface (c) When the Borel set is reduced to a single point on the surface, the asymptotic cone degenerates to the union of two planes}
 \label{fig::inout}
\end{figure}

We remark that the dimension of each fiber may change with $m$. This bundle is defined over the whole $M^N$, even at the points $m$ which are ``far" from ${\cal G}$. If $B(m,r)$ does not intersect the support of ${\cal G}$, then
${\cal C}_{B(m,r)}^{\text{par},{\cal G}}=T_m M^N$. This phenomenon is visualized in Figure \ref{fig::inout} where the cone \ref{subfconeleft} has its vertex (the black point) on the catenoid, and the cone \ref{subfconemiddle} has its vertex out of the catenoid.

\begin{figure}[H]
	\centering
	\begin{subfigure}[t]{0.48\textwidth}
                \includegraphics[width=\textwidth]{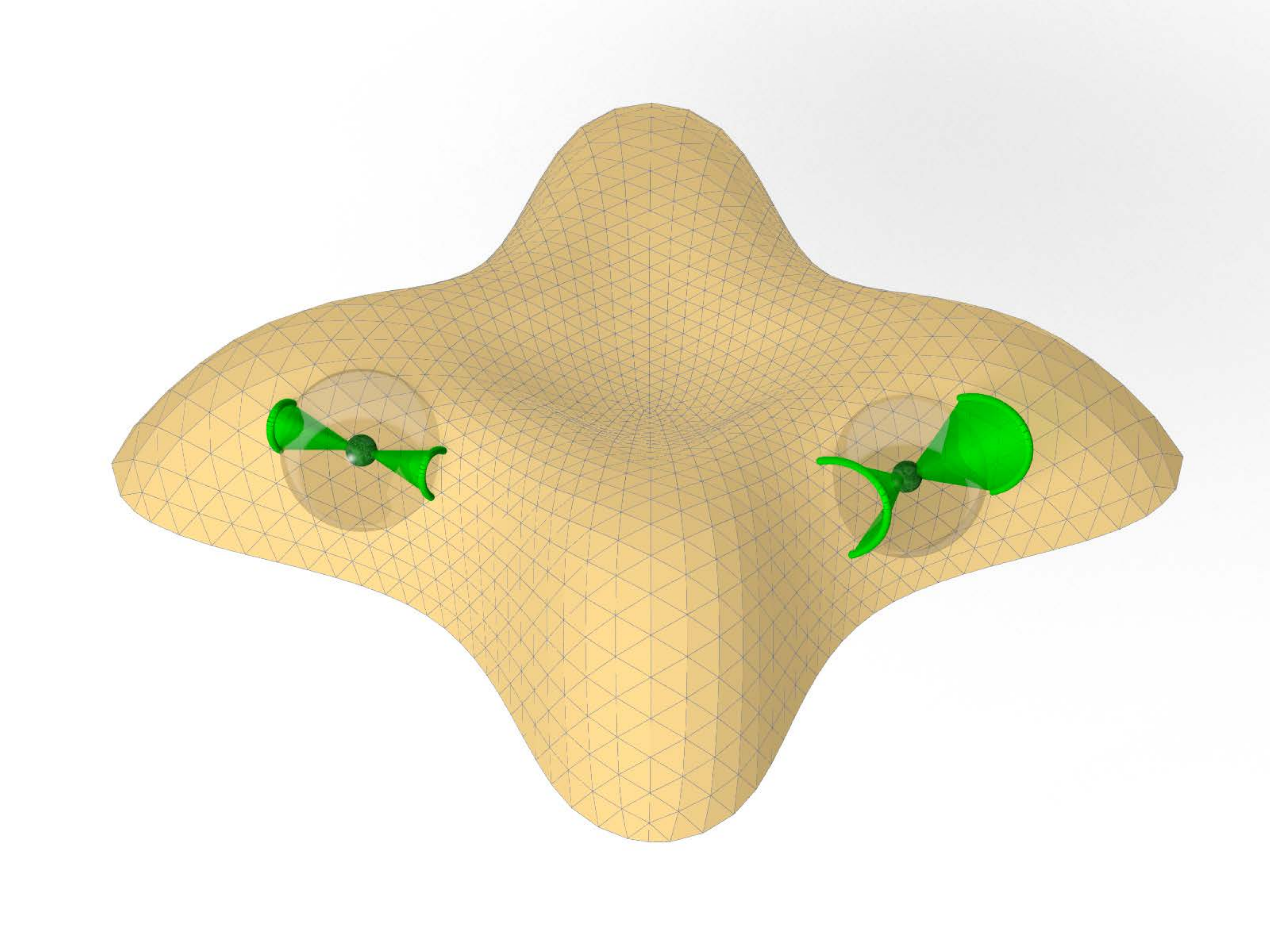}
      \end{subfigure}
	\begin{subfigure}[t]{0.48\textwidth}
        \includegraphics[width=\textwidth]{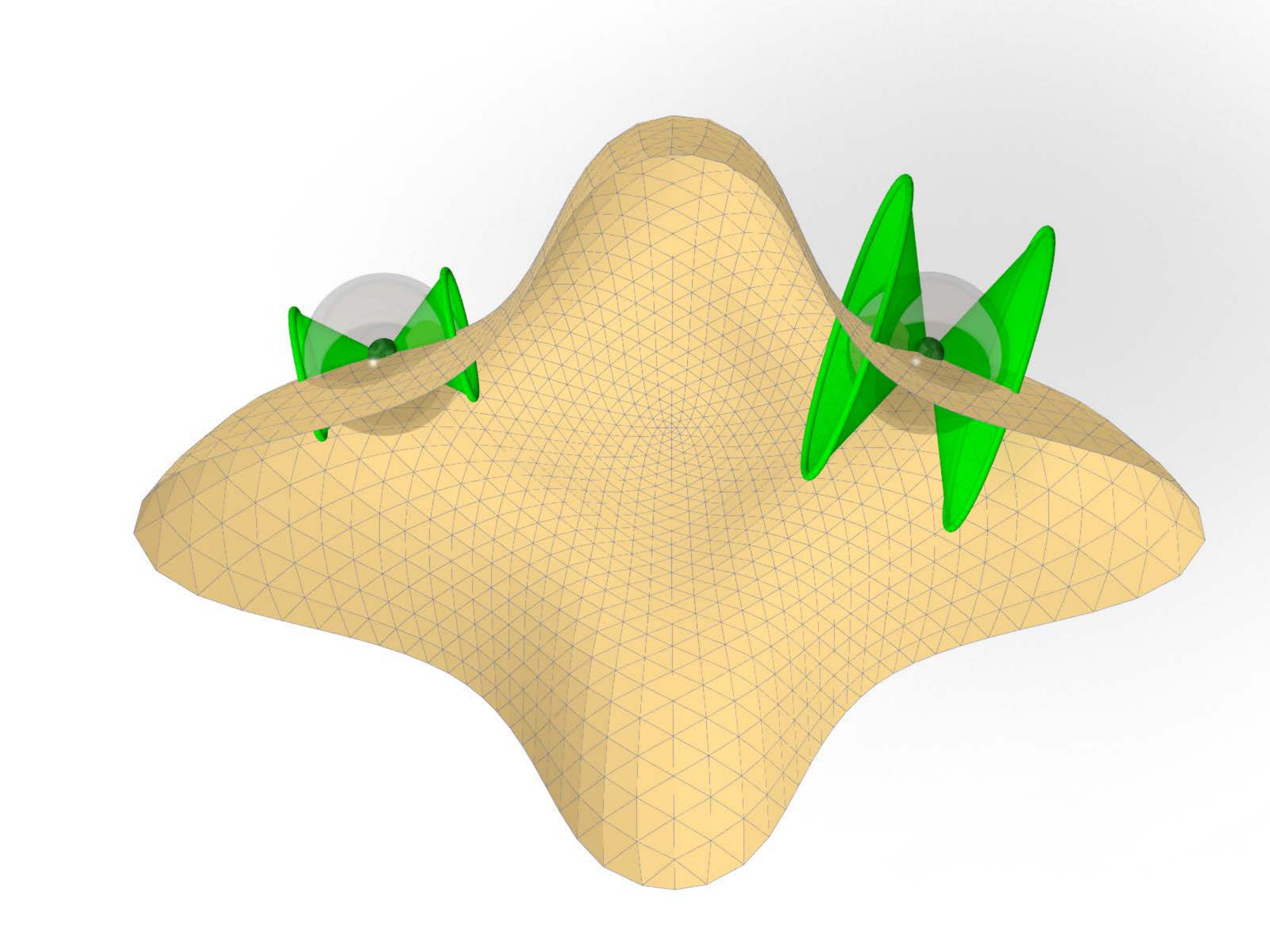}
    \end{subfigure}
    \caption{Asymptotic cones built on a non-smooth surface, here a triangulation: the top of the Lilium tower (designed by architect Zaha Hadid) in Warsaw, Poland.}
    \label{fig::liliumCone}
 \end{figure}

\begin{figure}[H]
	\centering
	\begin{subfigure}[t]{0.48\textwidth}
                \includegraphics[width=\textwidth]{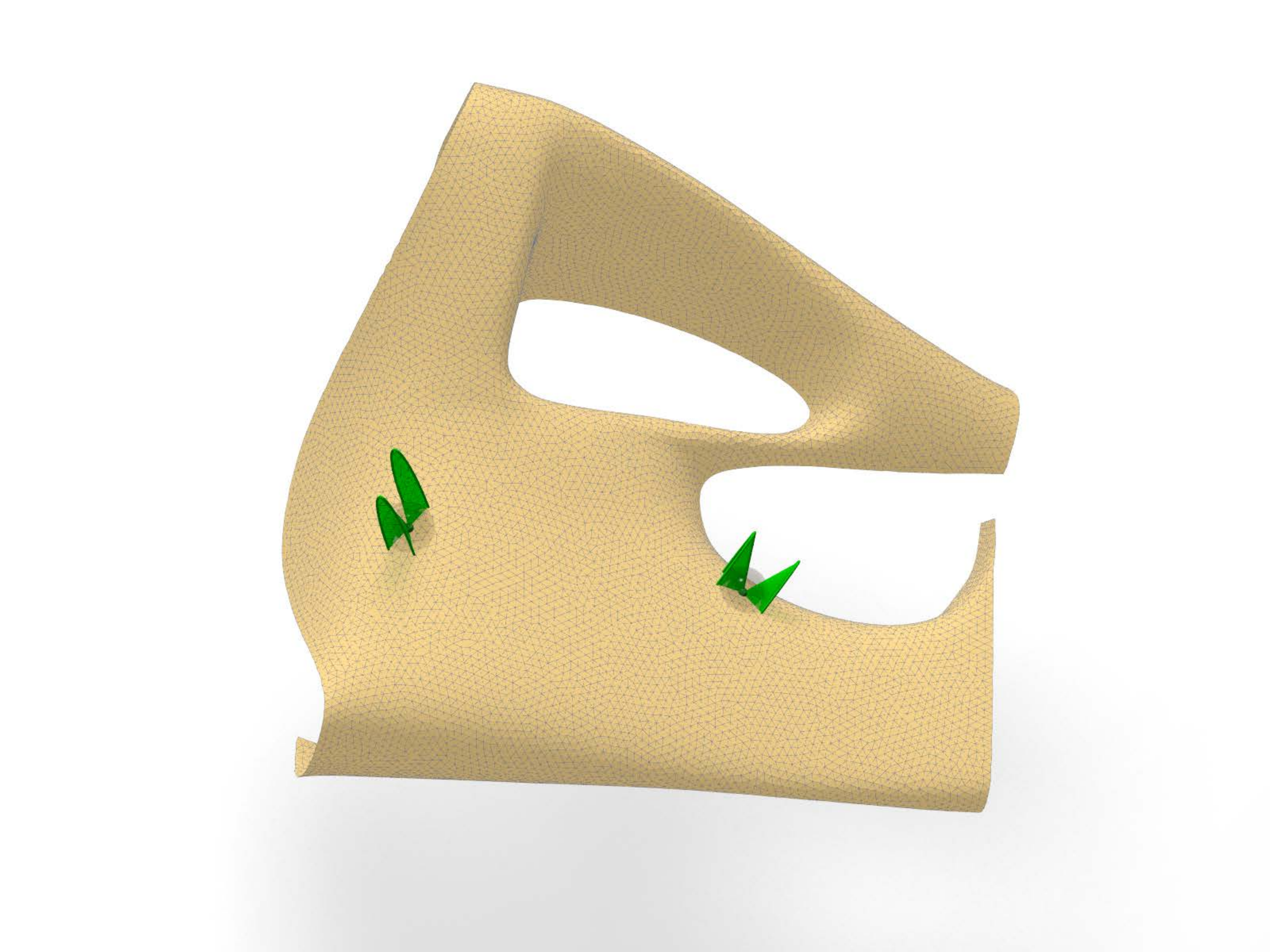}
      \end{subfigure}
	\begin{subfigure}[t]{0.48\textwidth}
        \includegraphics[width=\textwidth]{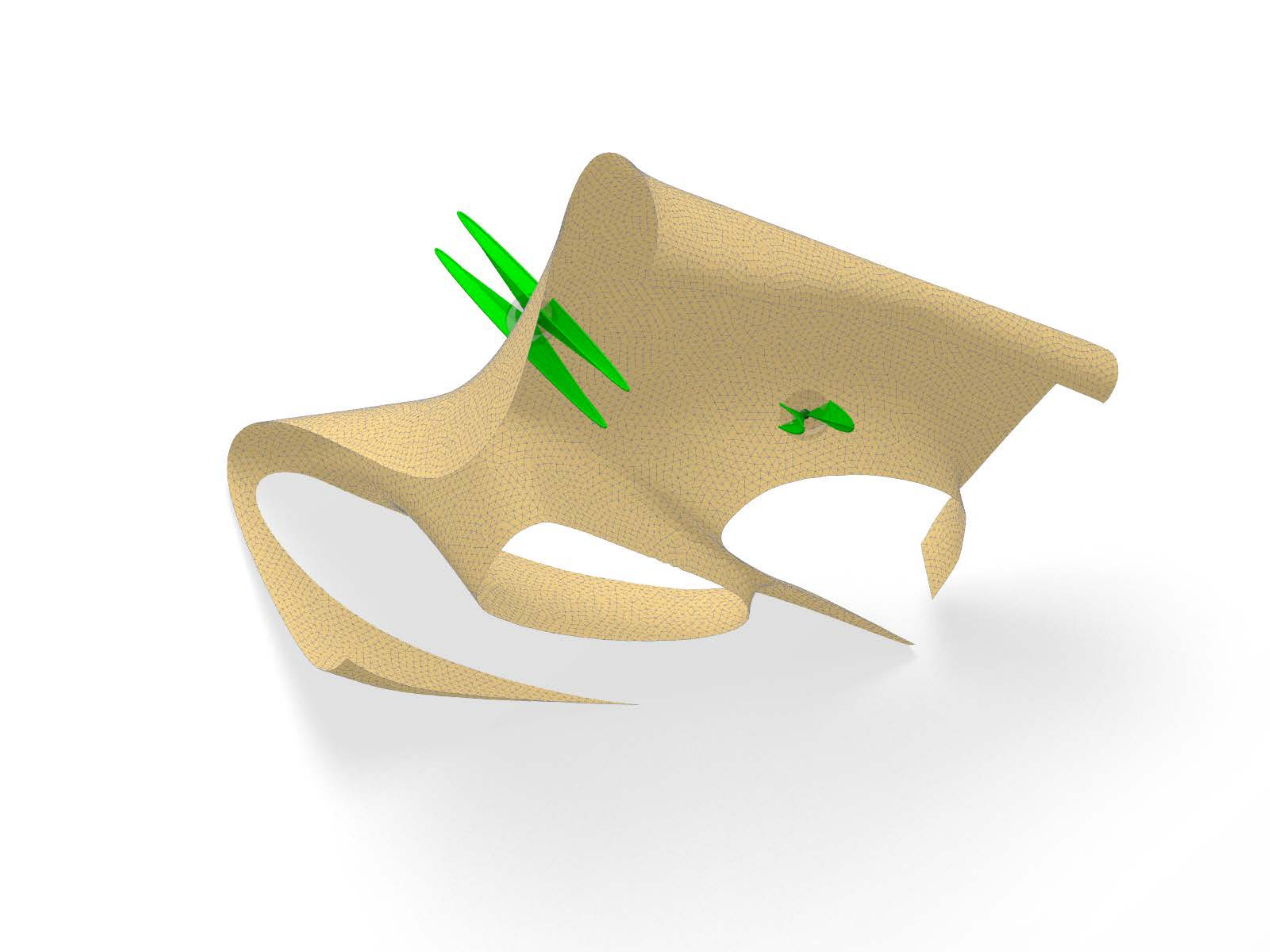}
    \end{subfigure}
    \caption{Asymptotic cones built on a non smooth-surface. Here another triangulation: the Heydar Aliyev Center (designed by architect Zaha Hadid) in Baku, Azerbaijan }
    \label{fig::bakuCone}
 \end{figure}

\section{Convergence and approximation results}
In the previous paragraphs, we gave our first
justification of the denomination of {\it asymptotic cones}. For surfaces $W^2$ in $\E^3$, this cone  reduces at each point $m$ of $W^2$ to  (the product of the normal line by) the standard asymptotic directions of $W^2$ at $m$.
We give now a second justification in terms of the convergence of sequences of polyhedra of $\E^N$. For simplicity, we restrict our study to (oriented) smooth hypersurfaces or polyhedra of $\E^N$ bounding a (compact) domain. We will show, in particular, that, if a sequence of domains ${\cal D}_k$ whose boundaries are polyhedra $P_{k}^{N-1}$, converges to a domain $D$ whose boundary is a smooth hypersurface $W^{N-1}$ (in a sense that will be clarified later), then for a large class of Borel subsets $B$, the sequence of asymptotic cones ${\cal C}_{B}^{{\cal D}_k}$ (resp.  ${\cal C}_{B}^{\text{par},{{\cal D}_k}}$) converges to ${\cal C}_{B}^{D}$ (resp.  ${\cal C}_{B}^{\text{par},D}$).
 We will use the following terminology \cite{morvan2008generalized}:
 \begin{itemize}
 \item
 The {\it fatness} $\Theta(P^{N-1})$ of a polyhedron $P^{N-1}$ is defined as follows: If $\sigma$ is an $l$-simplex, we begin to define the {\it size} $\epsilon(\sigma)$ of $\sigma$: It is the maximum over all edges $e$ of $\sigma$ of the length of $e$. Then, the {\it fatness} of $\sigma$ is the real number
$$\Theta(\sigma)= \min_{\mu \mbox{ {\footnotesize  $l$-simplex in }}  \sigma} \min_{j \in\{1,...,l\}} \frac{vol_j(\mu)}{\epsilon(\sigma)^j}.$$
Finally, the {\it fatness} of $P$ is the minimum of the fatness of its simplices. We denote by ${\bf F}_{\theta}$ the class of polyhedra in $\E^N$ with fatness greater or equal to $\theta$.

\item
  An $(N-1)$ dimensional polyhedron  $P^{N-1}$ in $\E^N$ is  {\it closely inscribed} in a smooth hypersurface $W^{N-1}$ if its vertices belong to $W^{N-1}$ and if the orthogonal projection of $P^{N-1}$ onto $W^{N-1}$ is a bijection.

\item
Let $P^{N-1}$ be an $(N-1)$ dimensional polyhedron (bounding a domain ${\cal D}$) closely inscribed in a smooth hypersurface $W^{N-1}$  in $\E^N$. The {\it angular deviation} $\alpha_m$ between $m \in P^{N-1}$ and  $pr_{W^{N-1}}(m)$ is the maximal angle between the normal $\xi_{pr_{W^{N-1}(m)}}$ of  $W^{N-1}$ and $(m,n_m)$, where $(m,n_m)$ belongs to the support of ${\bf N}({\cal D})$. If $B$ is any Borel subset of $\E^{N}$, we write
$$\alpha_B=\sup_{m \in B} \alpha_m.$$
\end{itemize}

One of the classical observations is that a ``good" fatness of a polyhedron and a ``small" Hausdorff distance of this polyhedron to a  smooth hypersurface in which it is closely inscribed imply that the angular deviation is ``small".
\subsection{A convergence result}
Let us state now our convergence theorem.
For any cone $\cal C$, we denote by   ${\cal UC}$ the basis of  ${\cal C}$; that is, ${\cal UC}={\cal C}\cap {\S}^{N-1}(0,1)$, where ${\S}^{N-1}(0,1)$ is the unit sphere centered at the origin.
We also use the distance $\tilde{d}$ defined on the class of  subsets of $\E^N$ by
\begin{equation}
\tilde{d}(A,B) = \inf_{x\in A, y \in B} d(x,y).
\end{equation}

\begin{theorem}\label{thm::triaconvergence}
Let $D$ be a (compact) domain  of $\E^N$ bounded by a smooth hypersurface $W^{N-1}$.  Let $({\cal D}_k)$ be a sequence of domains  of $\E^N$ bounded by polyhedra $(P_k^{N-1})$  closely inscribed in $W^{N-1}$ such that:
\begin{enumerate}
\item \label{HYP1}
The limit of $(P_k^{N-1})$ is $W^{N-1}$  for the Hausdorff distance.
\item \label{HYP2}
The fatness of $(P_k^{N-1})$ is uniformly bounded from below by a non-negative constant: there exists $\theta >0$ such that for all $k \in \N$, $P_k^{N-1} \in {\bf F}_{\theta}$.
\end{enumerate}
Let $B \in {\cal B}_{\E^N}$, such that for all $x\in \E^N$,  $|\Phi_{D}^{x}|(\partial B)=0$.   Then,
every sequence $(x_k \in {\cal C}_{B}^{{\cal D}_k})$ of unit vectors admits a subsequence (still denoted by $(x_k \in {\cal C}_{B}^{{\cal D}_k})$) that converges to a unit vector of  ${\cal C}_{B}^{D}$. In particular,
\begin{equation}
\lim_{k \to \infty}
\tilde{d}({\cal UC}_{B}^{{\rm par},{\cal D}_k},{\cal UC}^{{\rm par},D}_{B})=0.
\end{equation}
\end{theorem}

In the smooth case,  the assumption on the boundary of $B$ can be translated in terms of the second fundamental form of $W^{N-1}$. We say that the {\it normal curvature} at a point $p$ of a submanifold $V$ of $W^{N-1}$ is null if the second fundamental form $h_p$ of $W^{N-1}$ satisfies $h_p(u,u)=0$ for every $u$ tangent to $V$ at $p$. We get:

\begin{corollary}
Under assumptions \ref{HYP1} and \ref{HYP2} of Theorem \ref{thm::triaconvergence}, suppose that $B$ is a Borel set such that $W^{N-1} \cap \partial B$ is smooth and with null normal curvature. Then, the conclusion of Theorem \ref{thm::triaconvergence} holds.
\end{corollary}

To prove this corollary, we simply remark that under these assumptions,
\begin{equation}
\Phi_{D}^{x}(\partial B) = \int_{W^{N-1}\cap \partial B}h(pr_{TW^{N-1}}x,pr_{TW^{N-1}}x)=0.
\end{equation}

The rest of this section focuses on the proof of Theorem \ref{thm::triaconvergence}.

\subsubsection{Convergence of sequences of normal cycles}
We need to introduce the {\it flat norm} on  ${\cal D}_l(\E^N)$ as follows. The {\it mass} of an $l$-current $T$ is the real number
\begin{equation}
{\bf M}(T) = \sup \{T(\omega)\},
\end{equation}
where the supremum is taken over all $l$-differential forms with compact support such that $\sup_{m \in \E^N} |\omega_m| \leq 1$. The flat norm of an $l$-current $T$ is the real number
\begin{equation}
{\cal F}(T)= \inf \{{\bf M}(A) + {\bf M}(B)\},
\end{equation}
where the infimum is taken over all rectifiable $l$-currents $A$ and    $(l+1)$-currents $B$
such that $T= A+\partial B$.
Our main ingredient in our study of convergence and approximation of the asymptotic cones is the following result, which is a simple reformulation of Theorem $67$ of \cite{morvan2008generalized} (page $200$) for polyhedra:

\begin{theorem}\label{MAINTH}
If $P^{N-1}$ is a closed $(N-1)$ dimensional polyhedron bounding a domain $D$ and closely inscribed in a smooth closed hypersurface $W^{N-1}$ of $\E^N$ bounding a domain ${\cal D}$, then for any Borel subset $B$ of $P^{N-1}$,
\begin{equation}
\begin{split}
{\cal F}\big({\bf N}({\cal D})_{|T_B\E^N}& - {\bf N}(D)_{|T_{pr_{W^{N-1}}B}\E^N}\big)\\
\leq & K(\delta_B + \alpha_B){\bf M}({\bf N}({\cal D})_{|T_B\E^N}),
\end{split}
\end{equation}
where $\delta_B$ is the Hausdorff distance between $B$ and $pr_{W^{N-1}}(B)$ and $K$ is a constant depending on the norm of the second fundamental form of $W^{N-1}$.
\end{theorem}

The following proposition  can be deduced from Theorem \ref{MAINTH} in a slightly different version,
see also \cite{fu1993convergence}.

\begin{proposition}\label{FLATN}
Under the assumptions  \ref{HYP1} and \ref{HYP2} of Theorem \ref{thm::triaconvergence},
\begin{enumerate}
\item   \label{LCDF1}
The masses  ${\bf M}({\big(\bf N}({\cal D}_k))\big)$ are uniformly bounded from above;
\item  \label{LCDF2}
The sequence  ${(\bf N}({\cal D}_k))$  converges to  ${\bf N}(D)$ for the flat norm.
\end{enumerate}
\end{proposition}

\subsubsection{Convergence of sequences of asymptotic measures}
Our framework is the space of  (signed) Radon measures on $\E^N$ with finite total variation,  endowed with the norm $||.||_1$ defined for every  $\mu$ (with finite total variation $|\mu|$)  by
$$||\mu||_1 = \int_{\E^N}d|\mu|.$$
It is well known that this space is the (topological) dual to the space ${\cal C}_c(\E^N)$ of continuous functions with compact support on $\E^N$, endowed with the norm $||.||$ defined by
$$||f||= \sup_{x \in \E^N} |f(x)|.$$
The space  of  (signed) Radon measures (with finite total variation) can also be endowed with the
topology of {\it weak convergence of measures}: a sequence of Radon measures $(\mu_k)$  on $\E^N$ (weakly) converges to $\mu$ if, for every continuous function $f$ with compact support on
  $\E^N$ ({\it resp.} $T\E^N$), $\mu_k(f)$ converges to $\mu (f)$.

\begin{proposition}\label{FLATN2}
Under the assumptions  \ref{HYP1} and \ref{HYP2} of Theorem \ref{thm::triaconvergence},
\begin{enumerate}
\item \label{LSS1}
For each vector $x \in\E^N$, the sequence of  measures $(\Phi_{{\cal D}_k}^x)$ converges to  $\Phi_{D}^x$ for the weak convergence of measures on $\E^N$;
\item \label{LSS2}
For each unit vector $x \in\E^N$, the sequence of  measures $(\Phi_{{\cal D}_k}^x)$ is $||.||_1$-bounded.
\end{enumerate}
\end{proposition}

\noindent {\bf Proof of Proposition \ref{FLATN2}:} Item \ref{LSS1} is a direct consequence of Proposition \ref{FLATN}, since the flat convergence of the sequence of normal cycles ${(\bf N}({\cal D}_k))$ implies the weak convergence of the measures $(\Phi_{{\cal D}_k}^x)$. Item \ref{LSS2} is an application of the Theorem of Banach-Steinhaus: Since $(\Phi_{{\cal D}_k}^x)$ converges to  $\Phi_{D}^x$ for the weak convergence of measures on $\E^N$, for each $f \in  {\cal C}_c(\E^N)$, $\sup_k|<\Phi_{{\cal D}_k}^x,f><+\infty$. The Theorem of Banach-Steinhaus then implies that
$||\Phi_{{\cal D}_k}^x||_1$ is uniformly bounded with respect to $k$; that is,  the sequence of  measures $(|\Phi_{{\cal D}_k}^x|)$ is $||.||_1$-bounded.


Let us now explain our assumption on the boundary of the Borel subset $B$ in Theorem \ref{thm::triaconvergence}.
 In general, the weak  convergence of the sequence  $(\Phi_{{\cal D}_k}^x)$ {\it does not} imply  the convergence of
  $(\Phi_{{\cal D}_k}^x(B))$  to  $(\Phi_{D}^x(B))$ for every Borel subset. Indeed,
 generically,  characteristic functions are not continuous. That is why we  restrict our study to a class of Borel subsets with suitable boundaries with respect to $D$, such that we can use the following
 general lemma (see \cite{MLA} Chapter $4$ for instance):

 \begin{lemma}\label{LS}
 Let $(\mu_k)$ be a sequence of  (signed)  Radon measures on $\E^N$ such that
 \begin{enumerate}
 \item
  $(\mu_k)$ converges to $\mu$ for the weak topology,
  \item
  the sequence  $(|\mu_k|)$ of total variation of $(\mu_k)$ converges to a Radon measure $\nu$,  for the weak topology.
  \end{enumerate}
 If the boundary  $\partial B$ of  $B \in {\cal B}_{\E^N}$ satisfies $\nu(\partial B)=0$, then
\begin{equation}
\lim_{k \to \infty} \mu_k(B)= \mu(B).
\end{equation}
 \end{lemma}

%

 Since the sequence $(|\Phi_{{\cal D}_k}^x|)$ is $||.||_1$-bounded, we can extract a subsequence that converges. From Proposition \ref{FLATN} and Lemma \ref{LS},
 we deduce:

 \begin{proposition}\label{LAAS}
 Under assumptions  \ref{HYP1} and \ref{HYP2} of Theorem \ref{thm::triaconvergence}, if $x\in \Xi(E^N)$ and $B \in {\cal B}_{\E^N}$ satisfy  $|\Phi_{D}^{x}|(\partial B)=0$, then
 \begin{equation}
 \lim_{k \to \infty} \Phi_{{\cal D}_k}^x(B)=\Phi_{D}^x(B).
\end{equation}
 \end{proposition}

The last step of the proof of Theorem \ref{thm::triaconvergence} is to relate the behavior of  the sequence $(\Phi_{{\cal D}_k}^x)$ for any $x$ to the behavior of their associated quadratic cones.  For a fixed  $B \in {\cal B}_{\E^N}$, we will study  the quadratic forms
$x \to \Phi_{{\cal D}_k}^x(B)$
and
$x \to \Phi_{D}^x(B)$ introduced in Section \ref{ASMASC}.
We use the norm of uniform convergence on the space of quadratic forms: A sequence of quadratic forms $(q_k)$ defined on $\E^N$  {\it converges} to a quadratic form $q$ if $\sup_{||x||=1}|q_k(x)-q(x)|$ tends to $0$ when $k$ tends to infinity.

\begin{lemma}\label{lem::subsequence}
Let $(q_k)$ be a sequence of quadratic forms defined on $\E^N$, which converges to a quadratic form $q$.
Let $(x_k)$ be a sequence of unit vectors in $\E^N$, such that for each $k\in \N$, $x_k$ belongs to the isotropic cone ${\mathcal C}_k$ of $q_k$ ({\rm i.e.} $q_k(x_k)=0$). Then there exists a subsequence of $(x_k)$ that converges to a unit vector $x$ belonging to the isotropic cone ${\cal C}$ of $q$.
\end{lemma}

\noindent {\bf Proof of Lemma \ref{lem::subsequence}}:
We have
\begin{equation}\label{MKE}
\lim_{k \to  \infty}\sup_{||z||=1}|q(z) - q_k(z)|=0.
\end{equation}
Suppose that $(x_k)$ is a sequence of unit vectors such that for all $k$,  $x_k \in {\cal C}_k$.
Then, by the compacity of the unit sphere, there exists a subsequence (that we still denote by $(x_k)$), that converges to a unit vector $x$. From \ref{MKE}, we deduce
$$\lim_{k \to  \infty}|q(x_k) - q_k(x_k)|= \lim_{k \to  \infty}|q(x_k)|
=|q(x)|=0,$$
which means that $x \in {\cal C}$.


We remark that under the assumption of Lemma \ref{lem::subsequence}, we cannot claim that the sequence of Hausdorff distances between ${\cal UC}_k$ and ${\cal UC}$ tends to $0$ when $k$ tends to infinity, as shown in the following example.
Consider  the quadratic forms in $\E^3$ defined for every $k \in \N^*$ by
\begin{equation}
q_k(u,v,w)=\frac{1}{k}(u^2+v^2+w^2).
\end{equation}
We have for every $k\in \N^*$${\cal C}_k = \{ 0 \}$,
 ${\cal C} = \E^N$, ${\cal UC}_k=\emptyset$ and ${\cal UC}={\S}^{N-1}(0,1)$.


The proof of Theorem \ref{thm::triaconvergence} follows from Proposition \ref{FLATN}, Lemma \ref{LS},  Proposition \ref{LAAS} and Lemma \ref{lem::subsequence}.

\subsection{An approximation result}
In this section, for simplicity, we restrict our study to surfaces in $\E^3$.
We assume that $W^2$ (bounding $D$) and $P$ (bounding $\cal D$) are fixed, with $P$ being closely inscribed in $W^{2}$. We suppose (without any restriction) that $P$ is endowed with a triangulation, denoting by $t$ a generic triangle, and by  $r(t)$ its circum-radius. The following theorem compares the asymptotic cone of $\cal D$ over a Borel set $B$ composed of a union of triangles of $P$ and the asymptotic cone of $D$ over the orthogonal projection $pr_{W^{2}}(B)$ of $B$ on $W^{2}$.
If $C$  is the isotropic cone of a quadratic form $q$, we denote by
${\cal A}^{\epsilon}(UC)$ the set of unit vectors $x$ $\epsilon$-close to $C$; that is, $|q(x)| \leq  \epsilon$.

\begin{theorem}\label{APPRE}
Let $D$ be a domain  of $\E^3$ bounded by a smooth hypersurface $W^2$.  Let ${\cal D}$ be a domain  of $\E^3$ bounded by a polyhedron $P \in {\cal F}_{\theta}$, $\theta >0$  closely inscribed in $W^2$.
For any $\epsilon >0$, there exists $\eta >0$ such that if $\max \{r(t), t \in B\} \leq  \eta$,
then
$${\cal UC}^{{\rm par}, {\cal D}} \subset {\cal A}^{\epsilon}({\cal UC}^{{\rm par},D}_{{\rm pr}_{W^2}(B)}).$$
\end{theorem}
 In other words, roughly speaking, under the assumptions of Theorem  \ref{APPRE}, the asymptotic cones of ${\cal D}$ are close to the asymptotic cones of $D$. The proof  uses the results of Section $5$-$2$ of \cite{cohen2006second}. We summarize them in the following proposition.
\begin{proposition}\label{APPRE2}
Under the assumptions of Theorem \ref{APPRE}, for every unit vector $x$,
\begin{equation}\label{APPRE3}
|\Phi^{x}_{\cal D}(B)-\Phi^{x}_{D}({{\rm pr}_{W^2}(B)})| = {\bf K} \max \{r(t), t \in B\},
\end{equation}
where ${\bf K}$ is a constant depending on the area of $B$, the length of its boundary, the geometry of $W^{2}$, and  $\theta$.
\end{proposition}
If $x$ is a unit vector belonging to ${\cal UC}^{{\rm par}, {\cal D}}_B$, then \ref{APPRE3} implies that
$$|\Phi^{x}_{D}(B)| = {\bf K} \max \{r(t), t \in B\}.$$
Consequently, if the triangles of  $P$ have a sufficiently small circumradius, then
$$|\Phi^{x}_{D}({{\rm pr}_{W^2}(B)})| \leq \epsilon.$$
 The conclusion follows.

\section{Some experiments}

\subsection{Construction of asymptotic directions of a triangulation}
To mimic the smooth situation, it may be interesting to deduce asymptotic directions and asymptotic lines from the asymptotic cones defined on a singular surface. This construction may be achieved if one has a natural plane field on this surface. This is the case for triangulated surfaces, where each triangle spans a plane. Then, to build asymptotic directions at a point of a triangulated  surface $P^2$, associated to a Borel set $B$, one can  proceed as follows:
\begin{itemize}
\item
Consider a point $m$ on  $P^2$ (bounding $\cal D$) and build a Borel set $B$ around $m$; for instance, a ball centered at $m$ with a suitable radius.
\item
Build the asymptotic cone ${\cal C}_{B}^{{\rm par}, {\cal D}}$ whose vertex is $m$, associated with $B$.
\item
Build the intersection of  ${\cal C}_{B}^{{\rm par}, {\cal D}}$ with the plane spanned by the face of the triangle that contains $m$. The result is two directions, called the {\it asymptotic directions} of $P^2$ at $m$.
\item
When $m$ runs over the surface, build a cross field  (reduced to a point when the cone is reduced to a point). By integrating the cross field, we have asymptotic lines.
\end{itemize}

\begin{figure}[H]
	\centering
	\begin{subfigure}[t]{0.48\textwidth}
                \includegraphics[width=\textwidth]{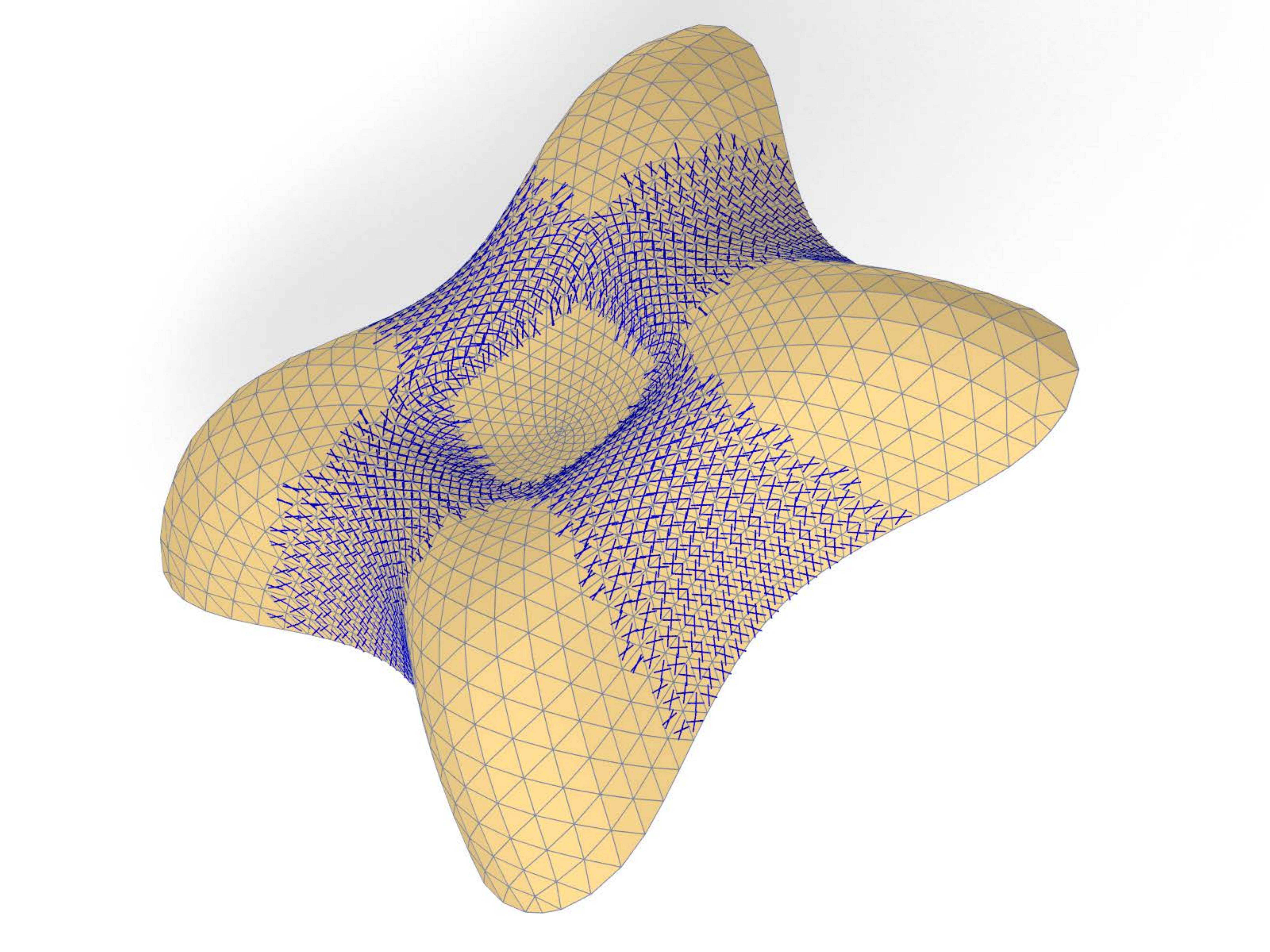}
      	 \subcaption{}
      \end{subfigure}
	\begin{subfigure}[t]{0.48\textwidth}
        \includegraphics[width=\textwidth]{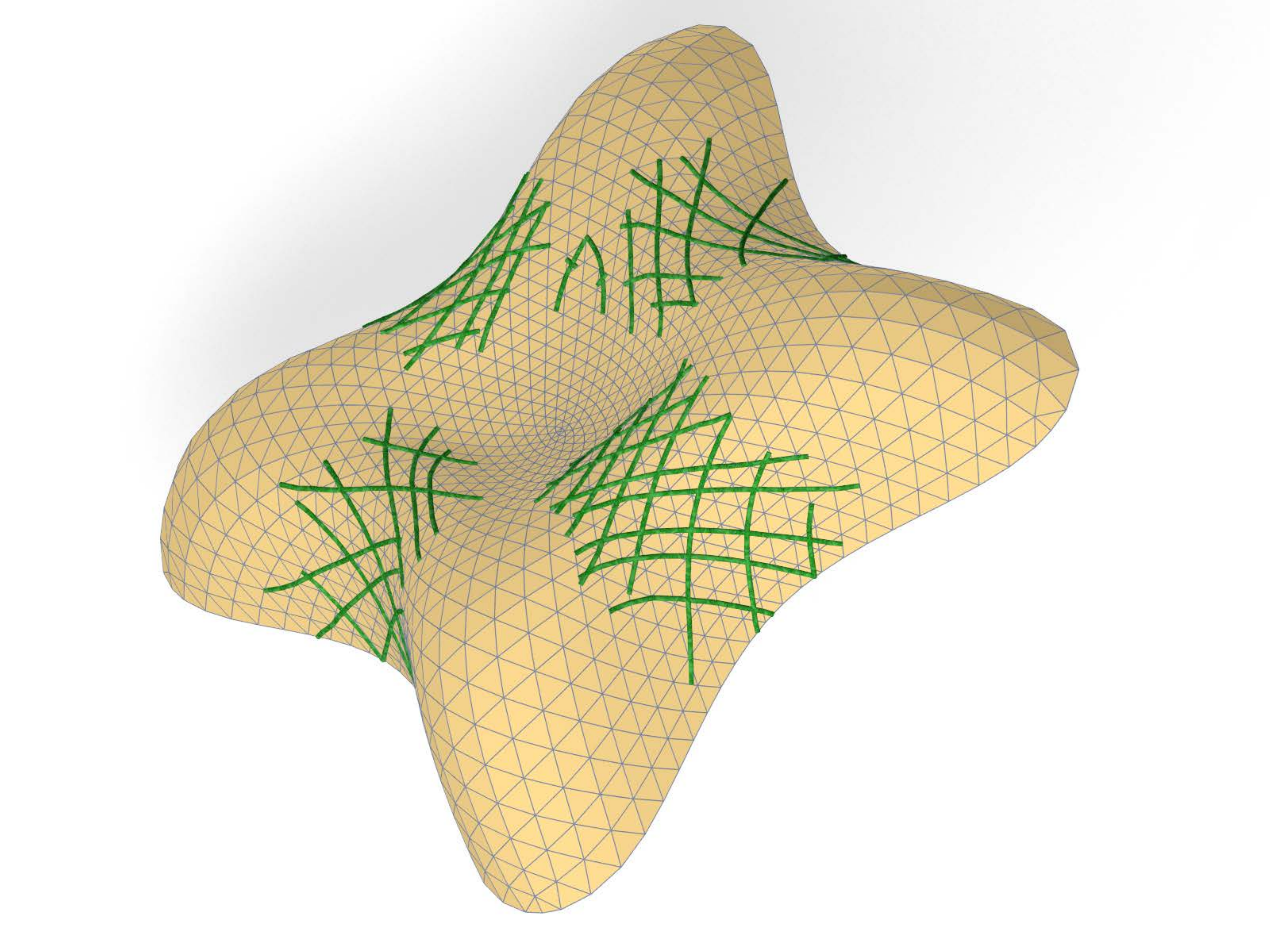}
        \subcaption{}
    \end{subfigure}
    \caption{(a) At each point where the asymptotic cone is not reduced to $\{0\}$, we build two asymptotic directions on the top of the Lilium tower (b) By integration we have asymptotic lines.}
    \label{fig::lilium}
 \end{figure}

\begin{figure}[H]
	\begin{subfigure}[t]{0.48\textwidth}
                \includegraphics[width=\textwidth]{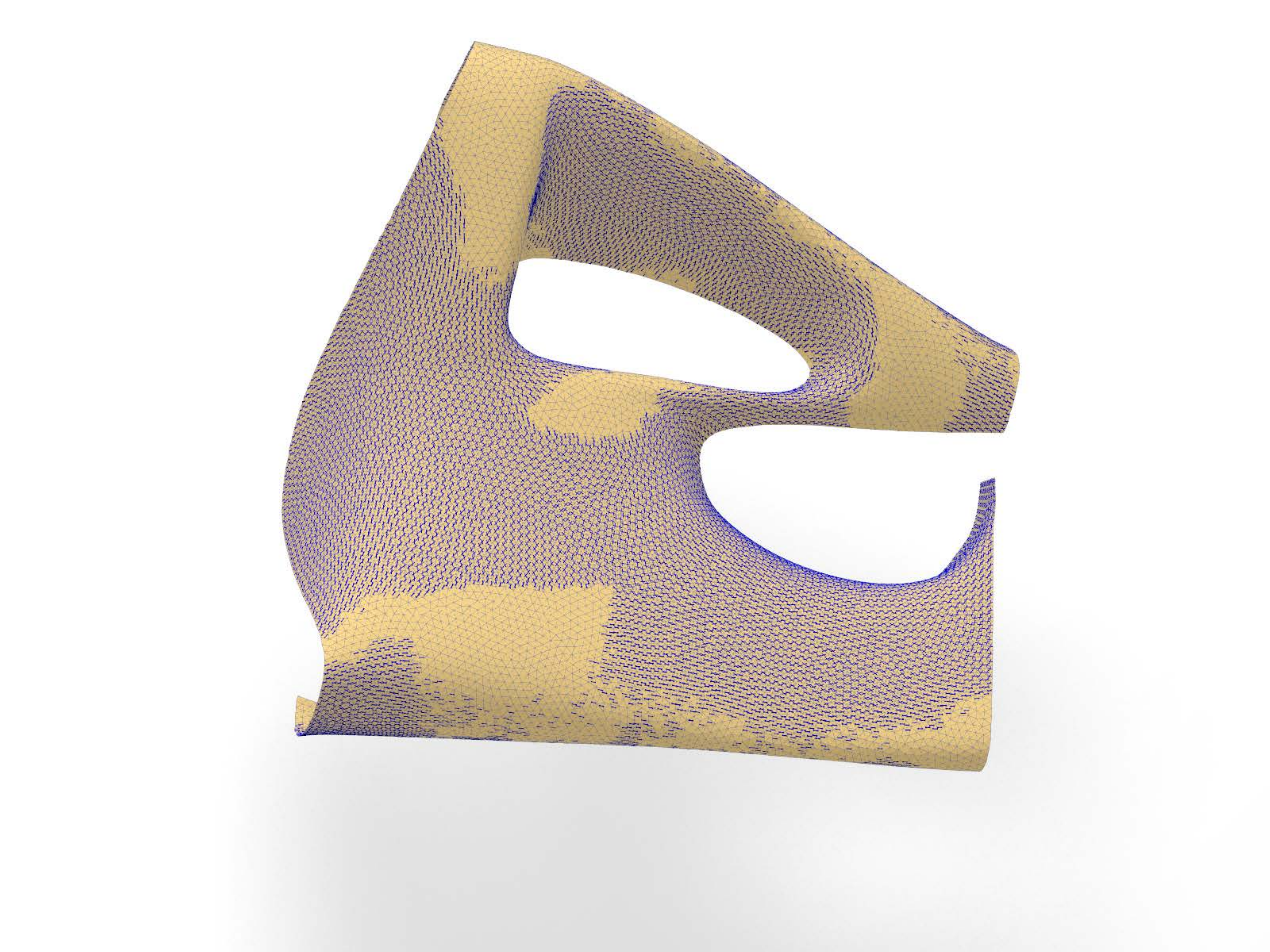}
      \subcaption{}
      \end{subfigure}
	\begin{subfigure}[t]{0.48\textwidth}
        \includegraphics[width=\textwidth]{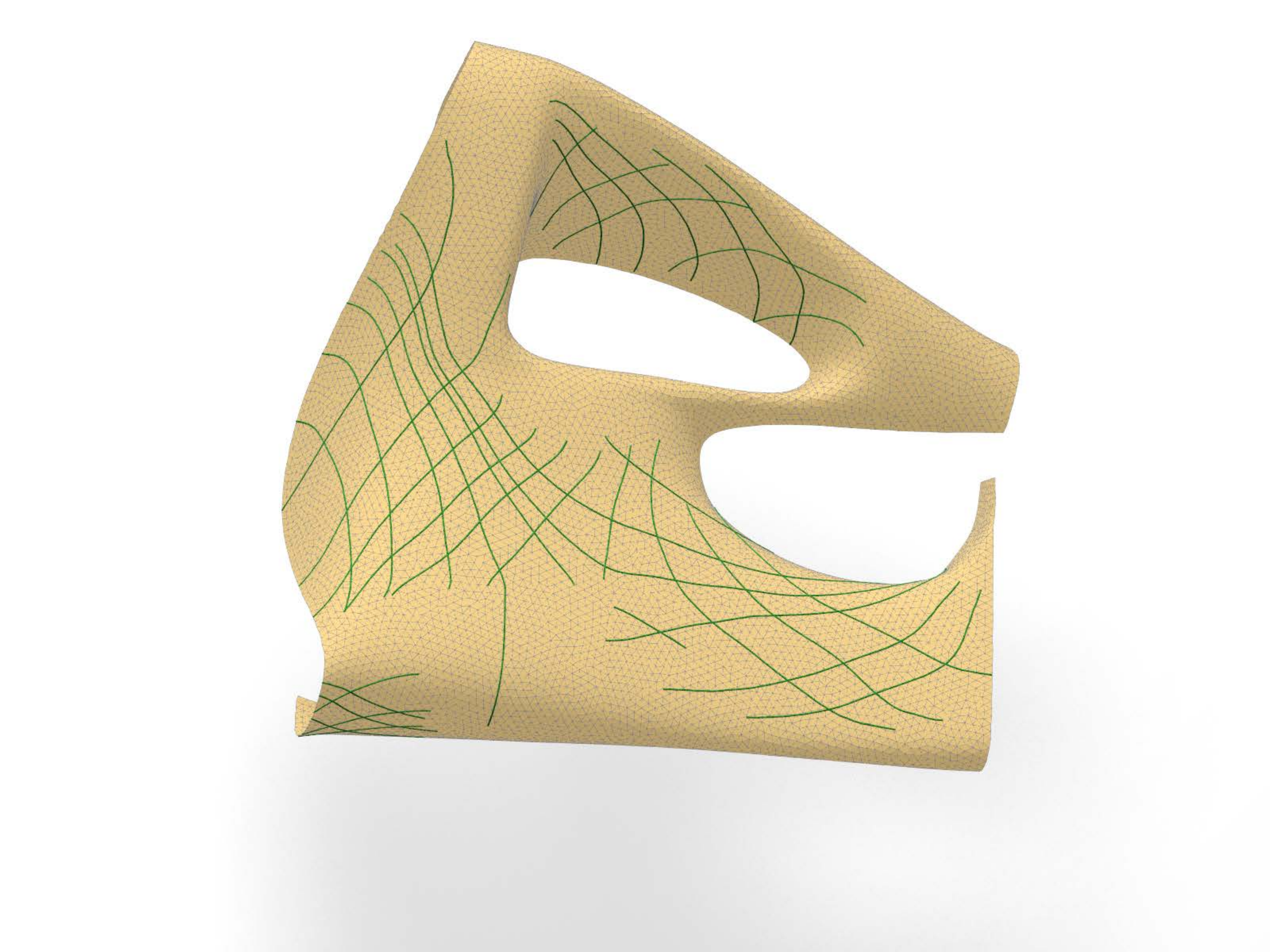}
        \subcaption{}
    \end{subfigure}
    \caption{(a) The cross field of asymptotic directions on the Heydar Aliyev Center in Baku (b) The  asymptotic lines obtained by integrating the asymptotic directions}
    \label{fig::baku}
 \end{figure}

\subsection{Approximation of the asymptotic lines of a smooth surface}
Using the construction of asymptotic cones, we can approximate the asymptotic directions ({\it resp.} lines) of a smooth surface $W^2$ in $\E^3$. We give an explicit example here.

\begin{figure}[H]
	\centering
	\begin{subfigure}[t]{0.48\textwidth}
                \includegraphics[width=\textwidth]{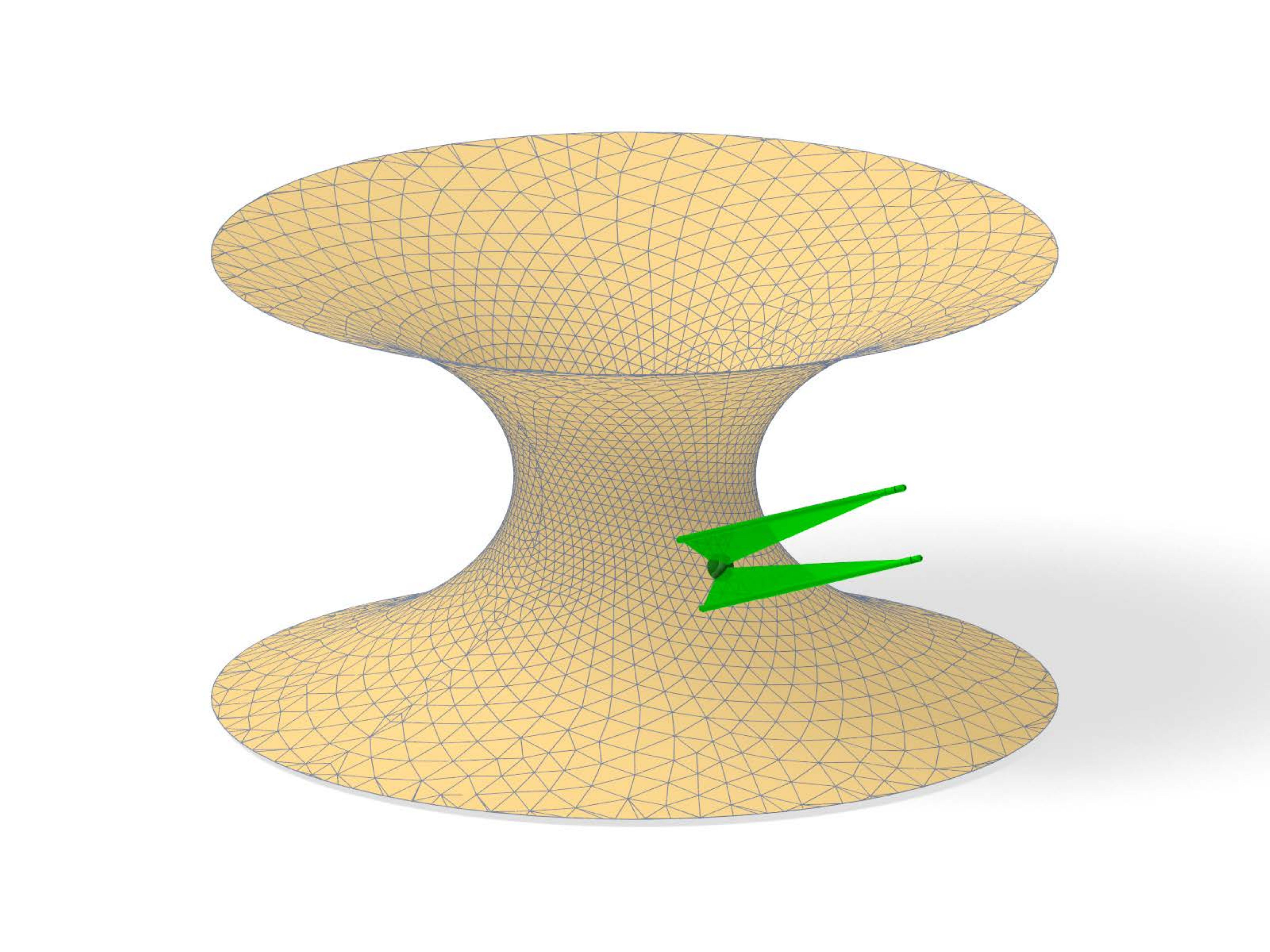}
                \subcaption{}\label{CATT1}
    \end{subfigure}
	\begin{subfigure}[t]{0.48\textwidth}
                 \includegraphics[width=\textwidth]{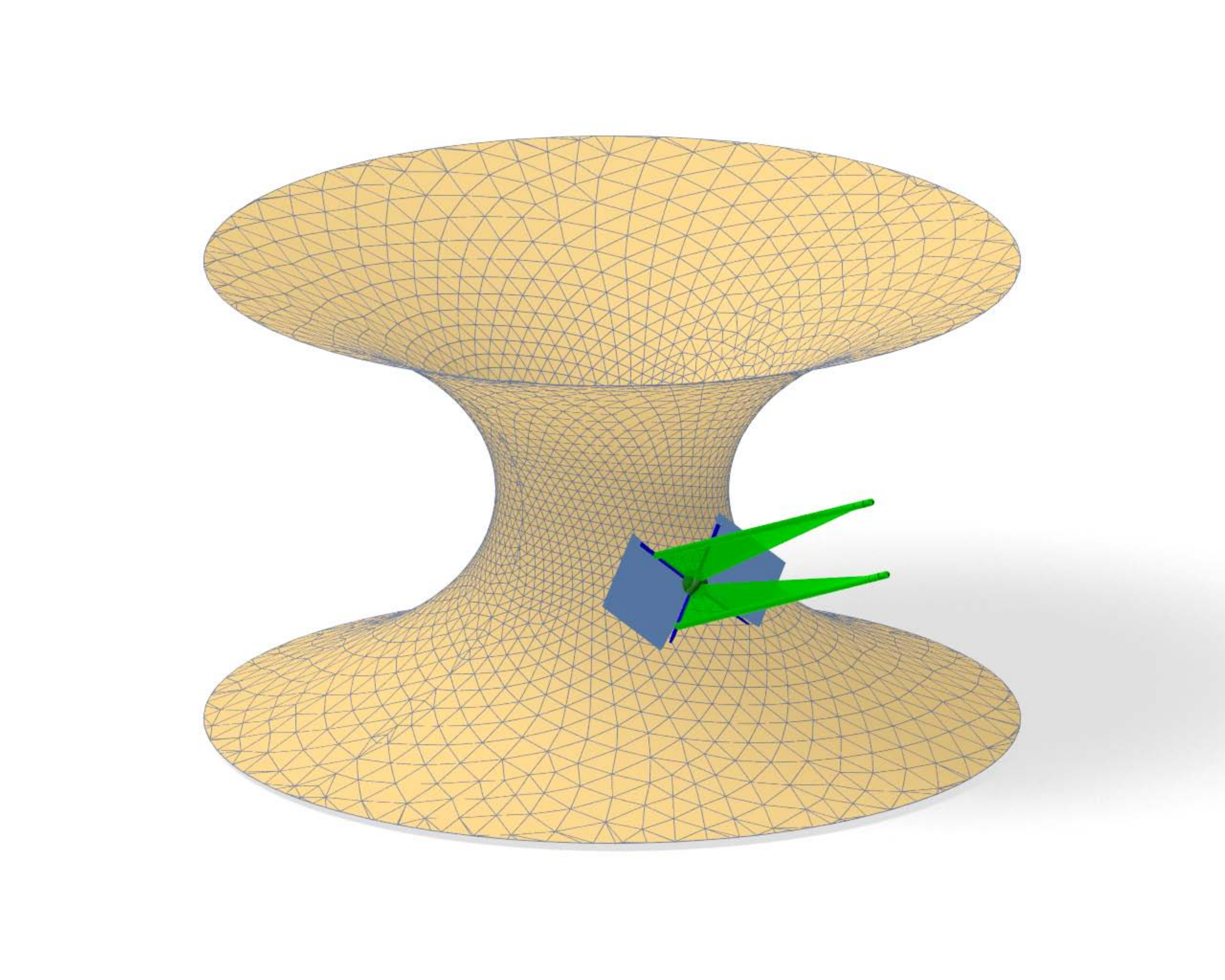}
                 \subcaption{}\label{CATT2}
    \end{subfigure}
	\begin{subfigure}[t]{0.48\textwidth}
                \includegraphics[width=\textwidth]{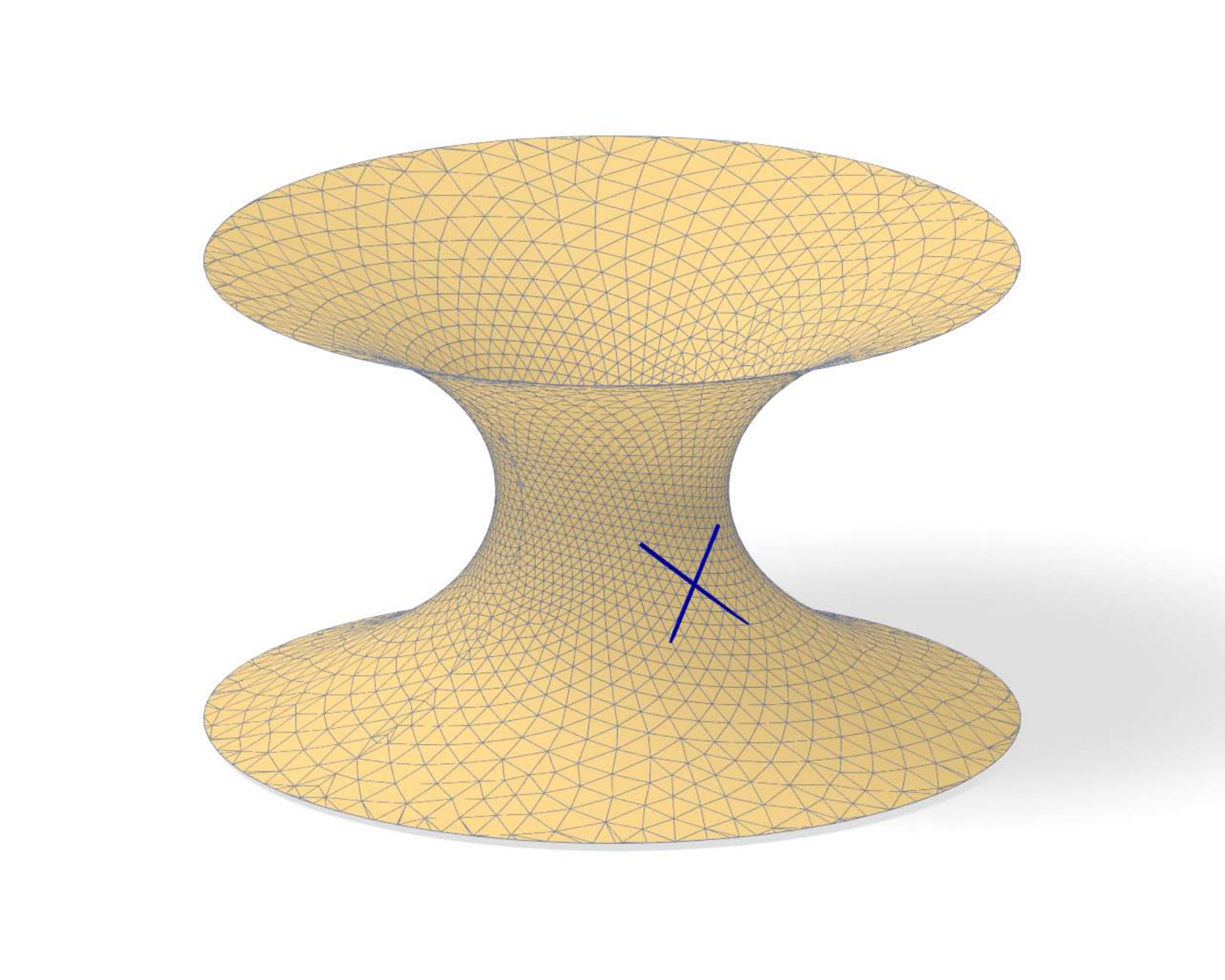}
                \subcaption{}\label{CATT3}
    \end{subfigure}
	\begin{subfigure}[t]{0.48\textwidth}
                \includegraphics[width=\textwidth]{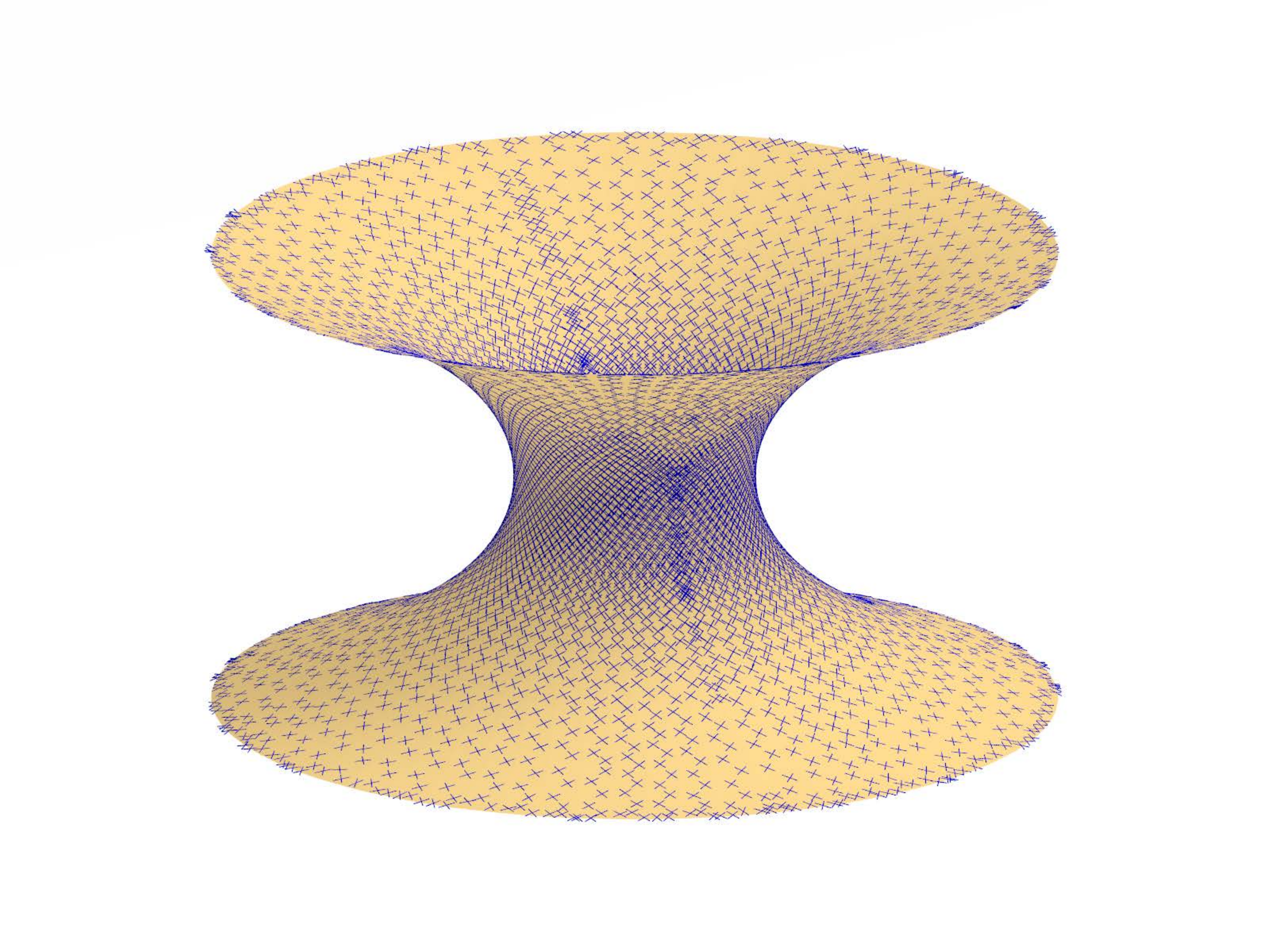}
                \subcaption{}\label{CATT4}
    \end{subfigure}
      \caption{(a) Step $1$: A green asymptotic cone built at a point $m$ of a triangulated surface  approximating a (smooth) catenoid (b) Step $2$: The blue plane is an approximation of the tangent plane of the (smooth) catenoid at $m$ (c) Step $3$: The intersection of the green asymptotic cone and the blue  plane gives two lines intersecting at $m$ (d) Step $4$: When $m$ runs over the triangulation, we obtain a  cross field  that approximates the field of asymptotic directions of the (smooth) catenoid}
\end{figure}

\begin{figure}[H]
\centering
	\begin{subfigure}[t]{0.48\textwidth}
    \includegraphics[width=\textwidth]{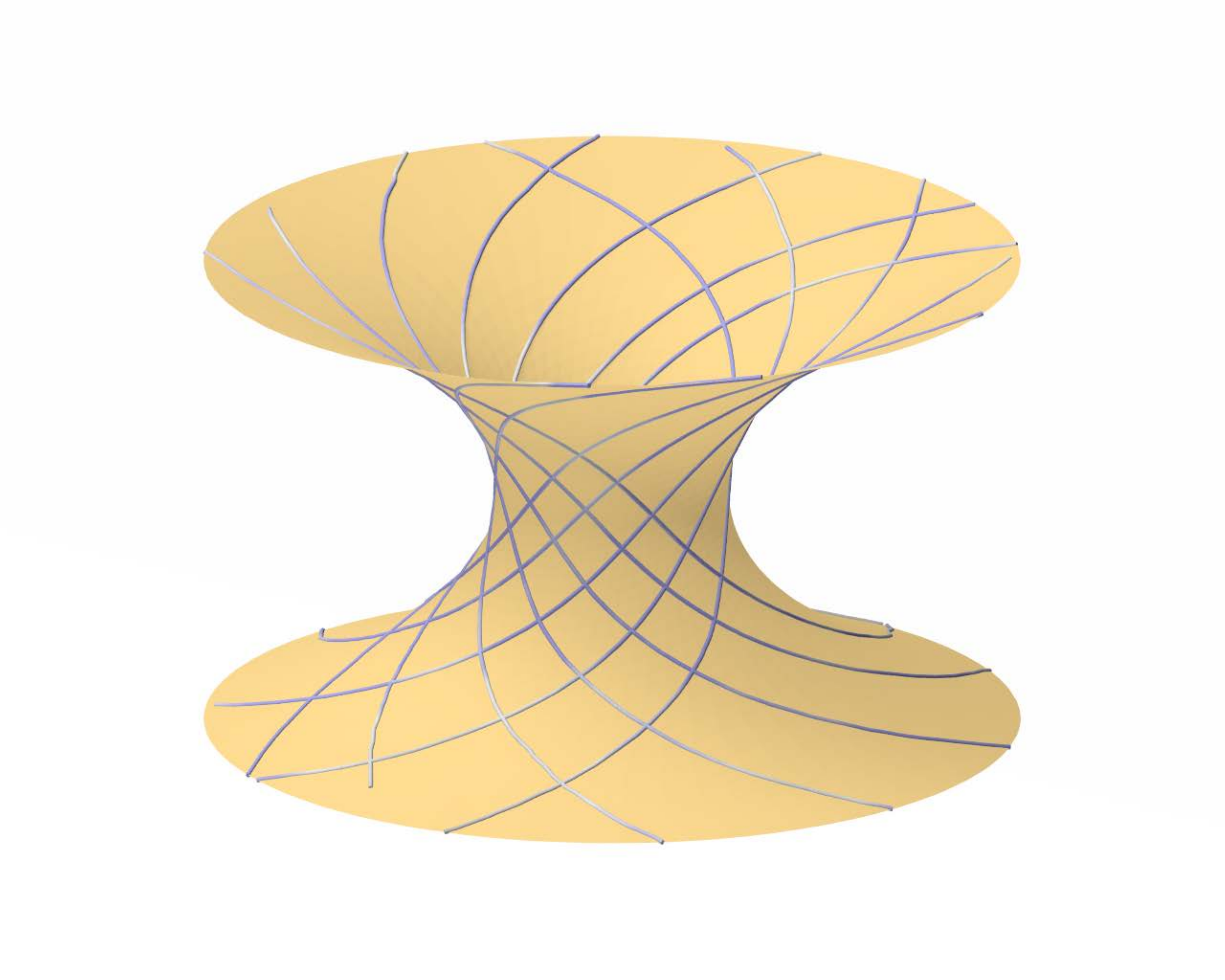}
     \subcaption{}\label{CATT5}
    \end{subfigure}
    \begin{subfigure}[t]{0.48\textwidth}
	\includegraphics[width=\textwidth]{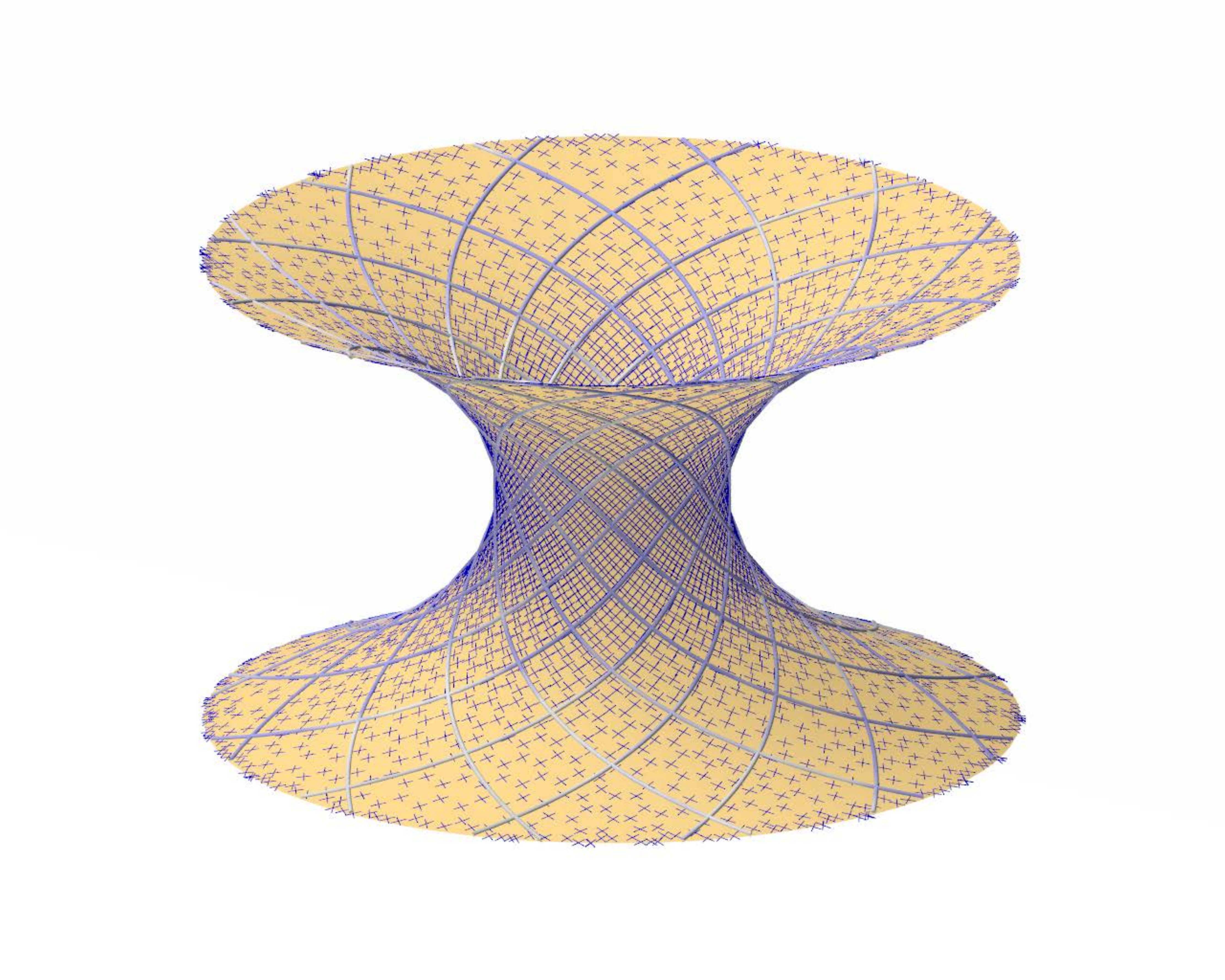}
     \subcaption{}\label{CATT6}
    \end{subfigure}
  \caption{(a) Integrating the cross field gives an approximation of the asymptotic lines of the catenoid (b) The cross field  is an approximation of the field of  asymptotic directions on the (smooth) catenoid. They can be compared with the blue lines, which are the asymptotic lines directly computed on the (smooth) catenoid}
\end{figure}

First of all, let us consider a portion of a  (smooth) catenoid   $W^2$, and a triangulation $P^2$  closely inscribed in  $W^2$, with a sufficiently dense set of vertices.  As shown in Figures \ref{CATT1},  \ref{CATT2},  \ref{CATT3},  \ref{CATT4}, the intersection of the tangent plane of $W^2$ at a point $m$ with the asymptotic cone of a ball centered at $m$ is reduced to two lines that are a discrete approximation of the asymptotic directions of $W^2$ at $m$. If necessary, we can approximate the tangent plane itself by the plane spanned by a face of the triangulation.
By integrating the directions field, we obtain  discrete asymptotic lines, as shown in Figures \ref{CATT5} and \ref{CATT6}. The reader can compare the asymptotic lines directly built on the smooth catenoid (these lines are orthogonal since the catenoid is a minimal surface) with the ``discrete ones" obtained by this procedure.

\subsection{Comparison of asymptotic lines}
We remark that the previous example gives (roughly speaking) ``very good" results
because the triangulated polyhedron approximating the catenoid suits it correctly, in the sense that the ``good" fatness of the triangles implies that the normal of any triangle $t$  is close to the normals of the orthogonal projection of $t$ on the catenoid.  To be more precise, we can estimate the error $er$  as follows: Let $W^2$ be a smooth surface approximated by a triangulated polyhedron $P^2$ closely inscribed on it. Let $m$ be a vertex of $P^2$, and  let us denote by $e_1$ and $e_2$ the asymptotic directions of $W^2$ at $m$, and by $\epsilon_1$, $\epsilon_2$ the approximated asymptotic directions at $m$. We define

 \begin{equation}\label{ERR}
 er_m = \inf \big(\frac{1}{2}(\angle(e_1, \epsilon_1) + \angle(e_2, \epsilon_2)), \frac{1}{2}(\angle(e_1, \epsilon_2) + \angle(e_2, \epsilon_1))\big),
 \end{equation}

where all the angles belong to $(0,\frac{\pi}{2})$.

\begin{figure}[H]
\centering
    \begin{subfigure}[t]{0.48\textwidth}
	\includegraphics[width=\textwidth]{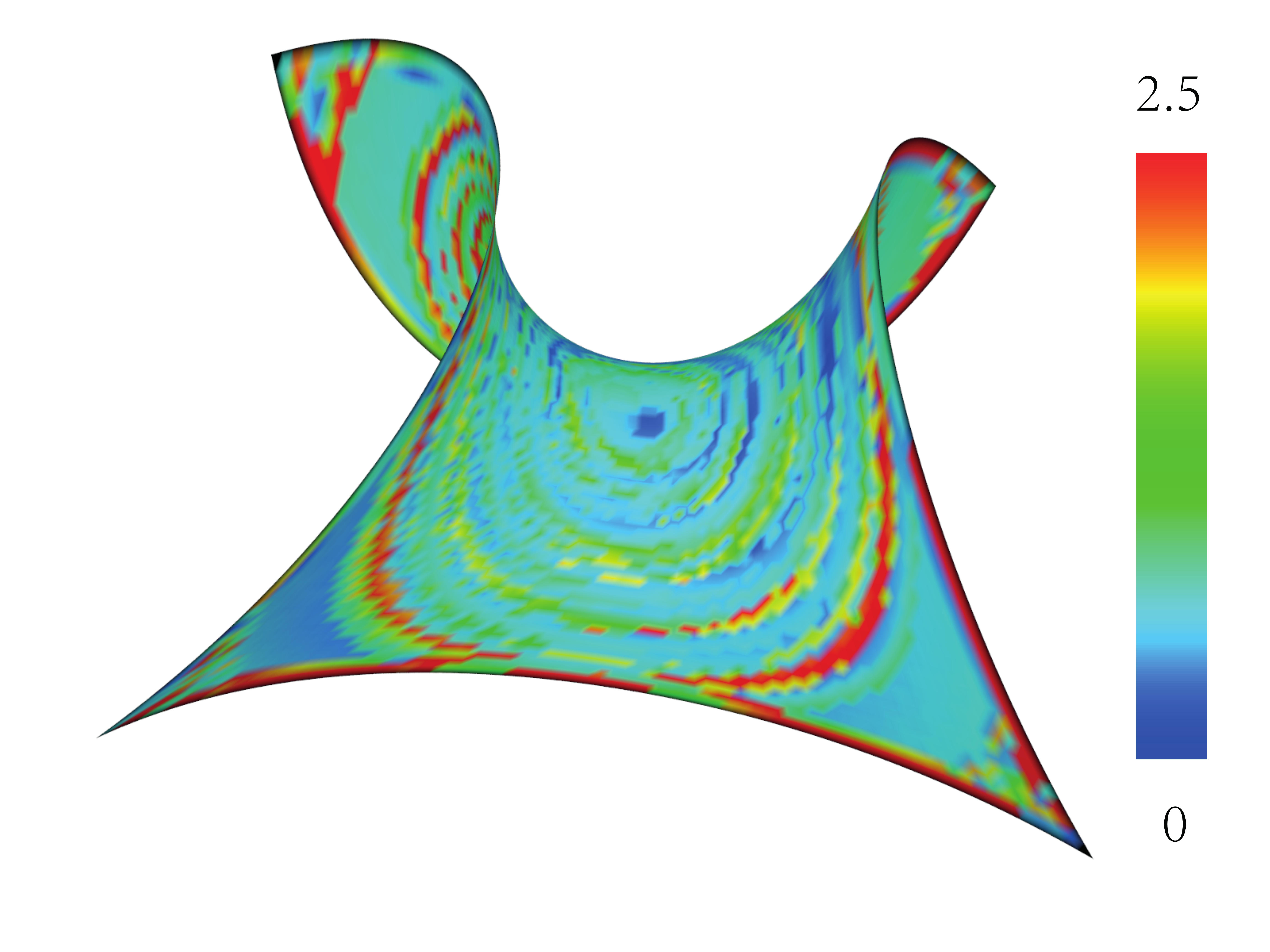}
	\subcaption{}
    \end{subfigure}
    \begin{subfigure}[t]{0.48\textwidth}
	\includegraphics[width=\textwidth]{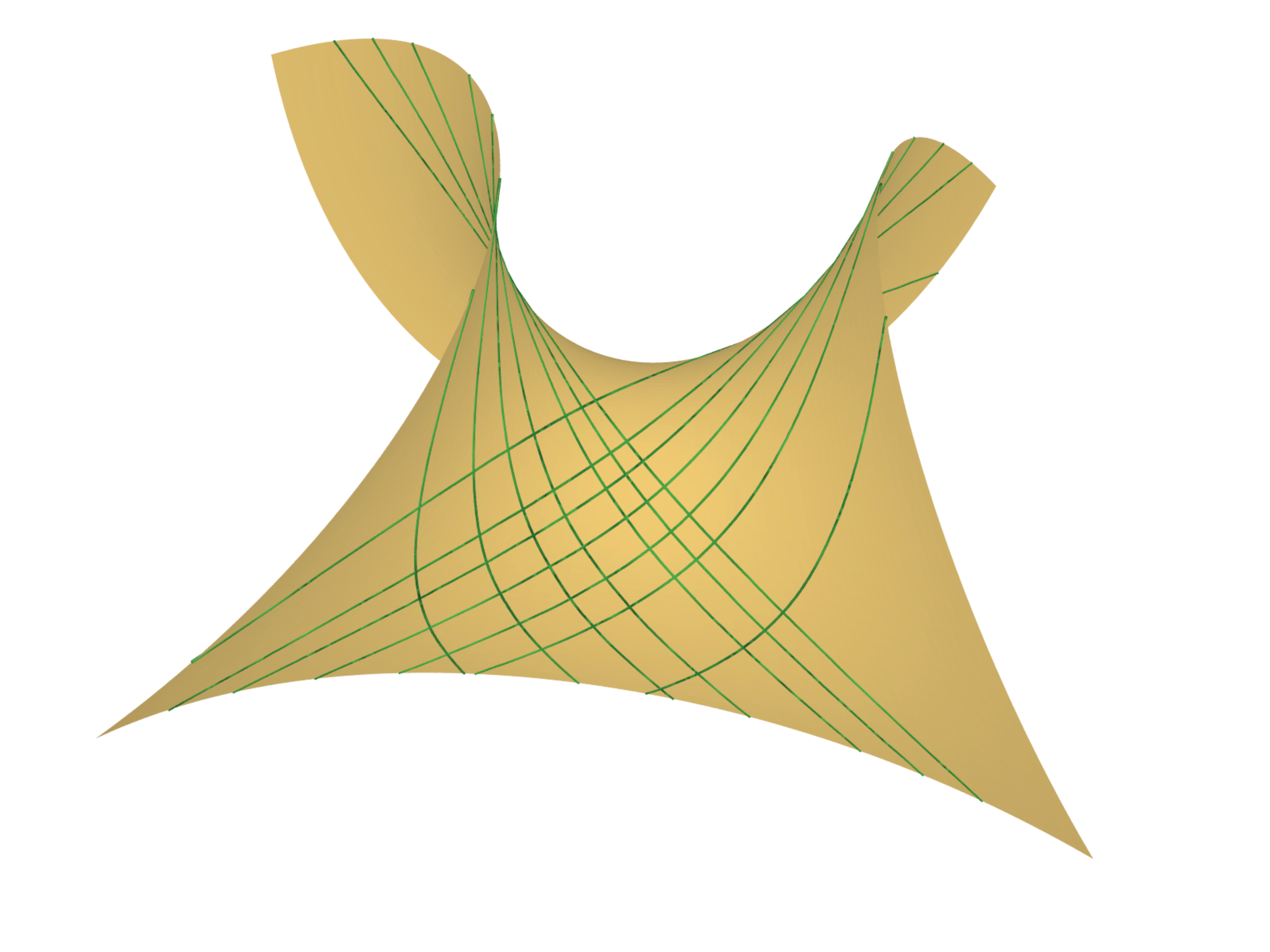}
	\subcaption{}
    \end{subfigure}
    \caption{(a) Comparison of discrete and smooth asymptotic lines on the Enneper surface (b) Some asymptotic lines computed on the polyhedron}\label{fig::enneper}
\end{figure}

For instance, let us consider the  portion of the Enneper surface shown in Figure \ref{fig::enneper}, and a triangulated polyhedron closely inscribed on it.
The error $er$ is always less than or equal to $5$ degrees.
\label{fig::enneper88}

\begin{figure}[H]
	\centering
    \begin{subfigure}[t]{0.48\textwidth}
	\includegraphics[width=\textwidth]{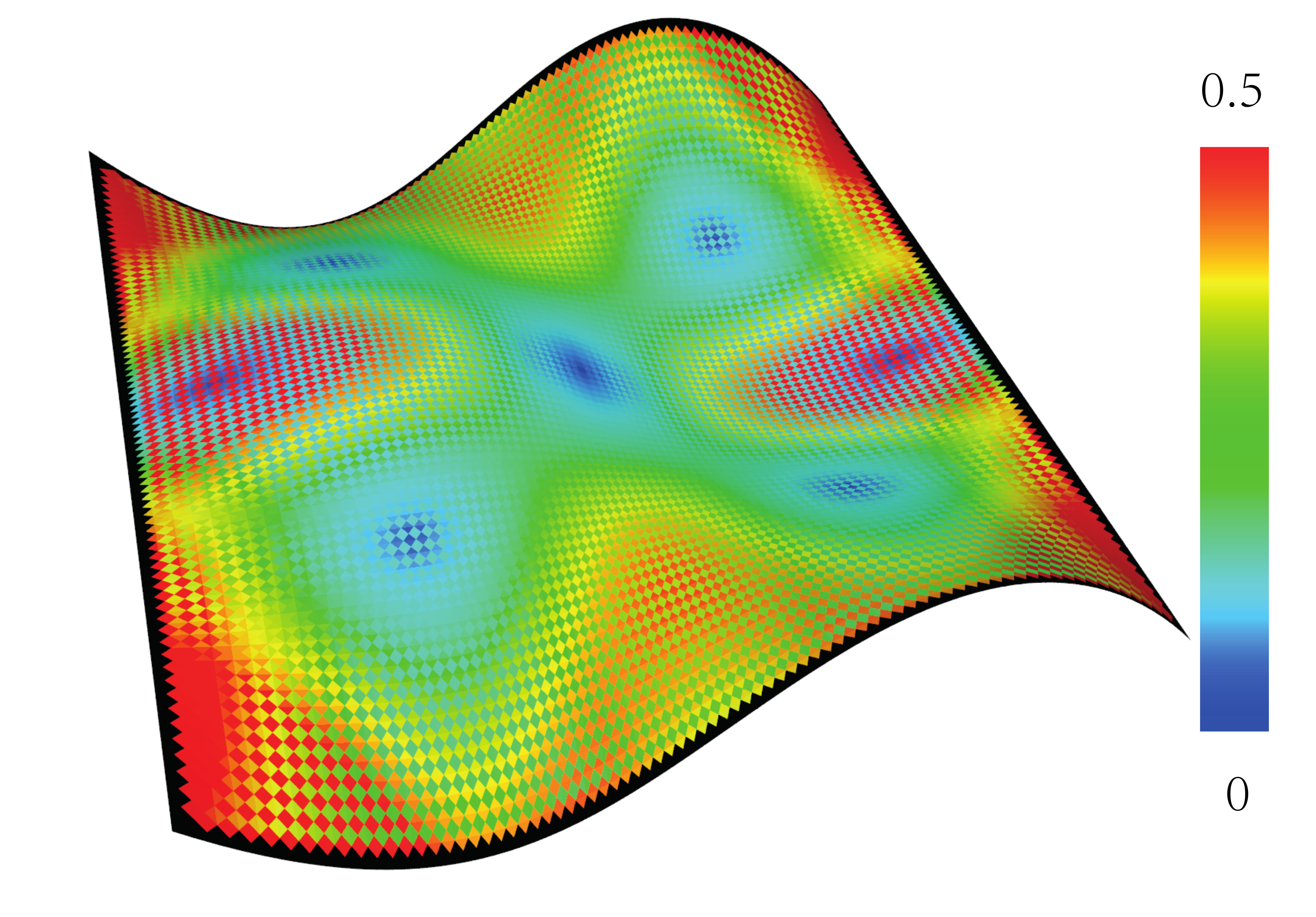}
	\subcaption{}
    \label{fig::splinesurfaceNormalsCoarse}
    \end{subfigure}
    \begin{subfigure}[t]{0.48\textwidth}
	\includegraphics[width=\textwidth]{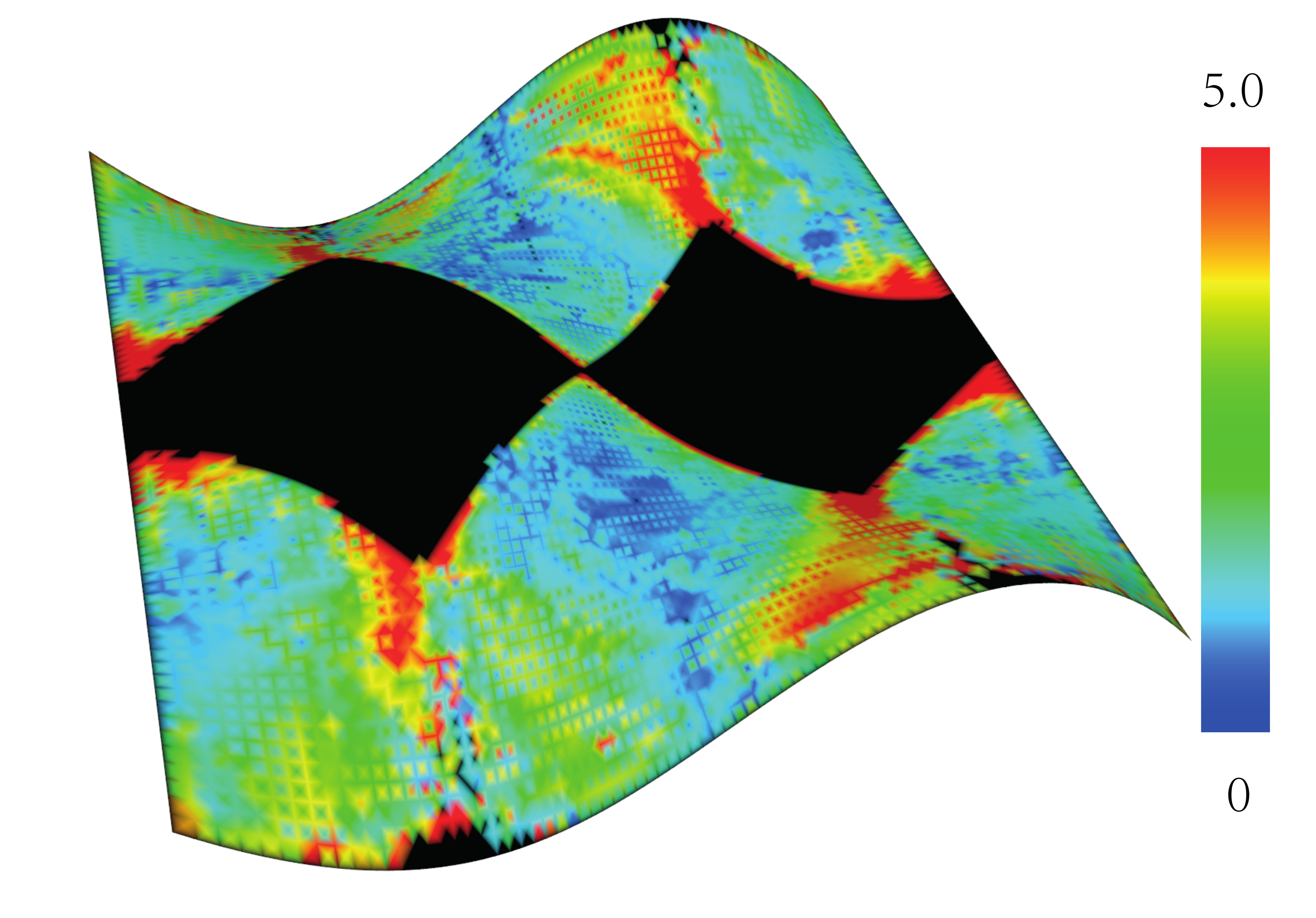}
	\subcaption{}
    \label{fig::splinesurfaceAngleDifference}
    \end{subfigure}
    \begin{subfigure}[t]{0.48\textwidth}
	\includegraphics[width=\textwidth]{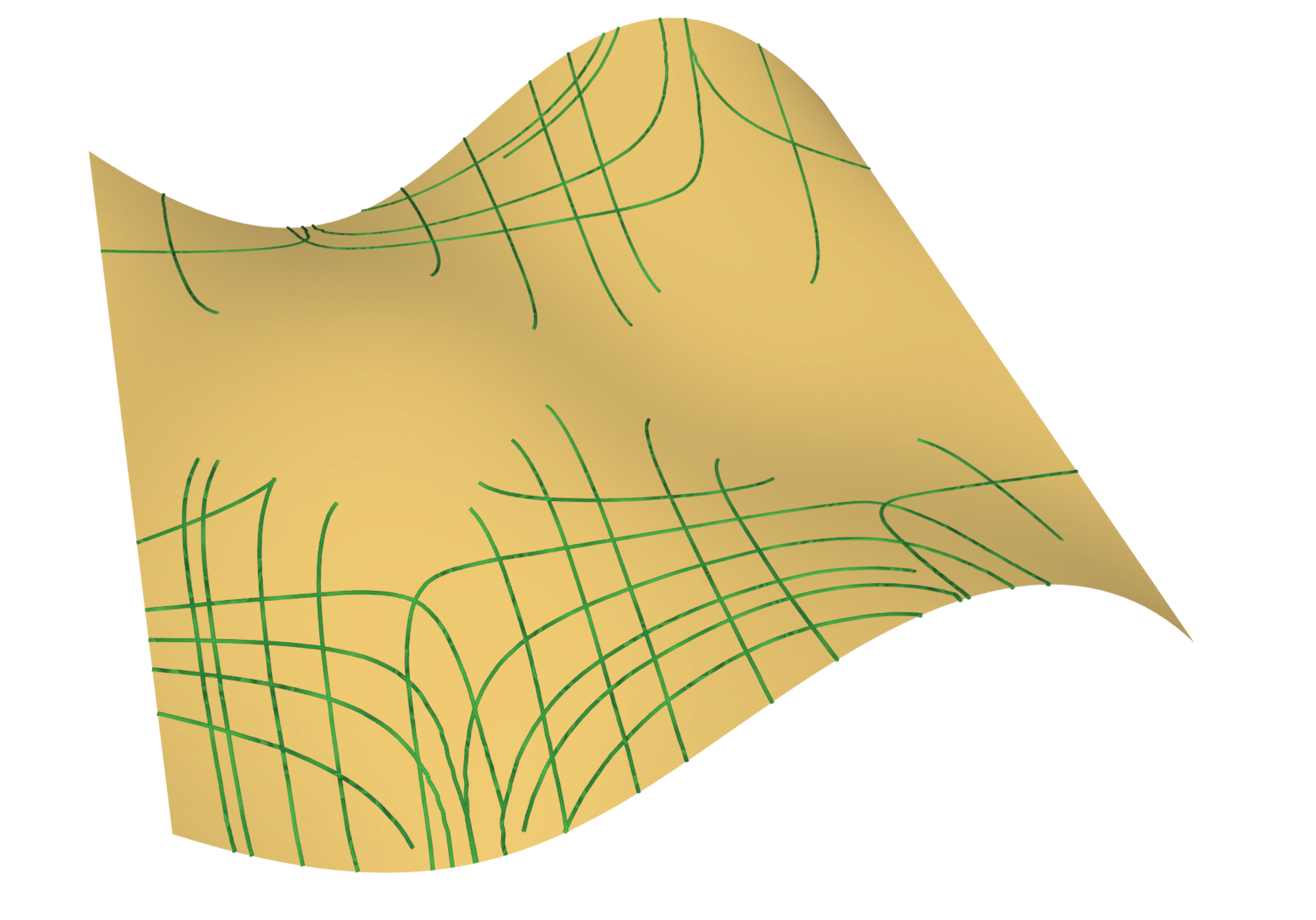}
	\subcaption{}
    \label{fig::splinesurfaceLines}
    \end{subfigure}
    \caption{(a)  Comparison of the normals computed on a $B$-spline surface and an approximating triangulation (b) Comparison of the asymptotic directions by computing $er_m$ at each vertex $m$  (the black points are convex points, where the asymptotic cone is reduced to $\{0\}$ (c) Some asymptotic lines computed on the triangulation }
\end{figure}

In the following example, we show that, in accordance with the theory, the error may be large even with a very thin triangulation closely inscribed in a smooth surface, if the angle between the tangent plane at a point and the corresponding triangle is too large (the same phenomenon appears when one compares  the area of a cylinder with the area of a Lantern of Schwarz inscribed on it, see \cite{morvan2008generalized} for instance). Here, we consider a (smooth) surface $W^2$ in $\R^3$ obtained as the graph of a tensor product $B$-spline function
$$f : G=[0,1] \times [0,1] \to \R,$$
of degree $3$ in each variable, defined by $5 \times 5$ control points over a regular grid of $80 \times 80$ points. Each square of $G$ is triangulated by taking both diagonals. We build  the corresponding (piecewise linear) triangulation $P^2$ inscribed on $W^2$.  Then, we compare, at each vertex $m$ of the triangulation, the normal of $W^2$ at $m$ and the average of the normals of the triangles incident to $m$. The error varies between $0$ and $0.5$ degree (see Figure \ref{fig::splinesurfaceAngleDifference}). On the other hand, we compute on the same triangulation the  asymptotic directions of $W^2$ and the asymptotic cones of $P^2$, from which we deduce discrete asymptotic directions by intersecting the cones with the tangent planes as before. Then, we compare at each point $m$ the error $er_m$ given in Figure \ref{fig::splinesurfaceNormalsCoarse}. The error at each vertex $m$ varies between $0$ and $5$ degrees, according to the behavior of the normal of the faces incident to $m$  with respect to the normal of $W^2$ at $m$.

\label{LJK7}


\subsection{Deformation of asymptotic lines of discrete surfaces}
In the following example, we produce a deformation of ``discrete" asymptotic lines as follows: We build a triangulation closely inscribed on a smooth surface  $W^2$ (here a Chen's surface). We then compute asymptotic lines by the previous process using balls of radius $R=3$ (the normalization is such that the average length of the edges is $1$). Then, we slightly modify the position of the vertices that can now be out of  $W^2$ (in other words, we create noisy data), without modifying the topology of the triangulation. With this new triangulation, we compute new asymptotic lines using the same process and using balls of the same radius ($R=3$) or a different radius ($R=6$).

 \begin{figure}[H]
	\centering
	\begin{subfigure}[t]{0.31\textwidth}
    \includegraphics[width=\textwidth]{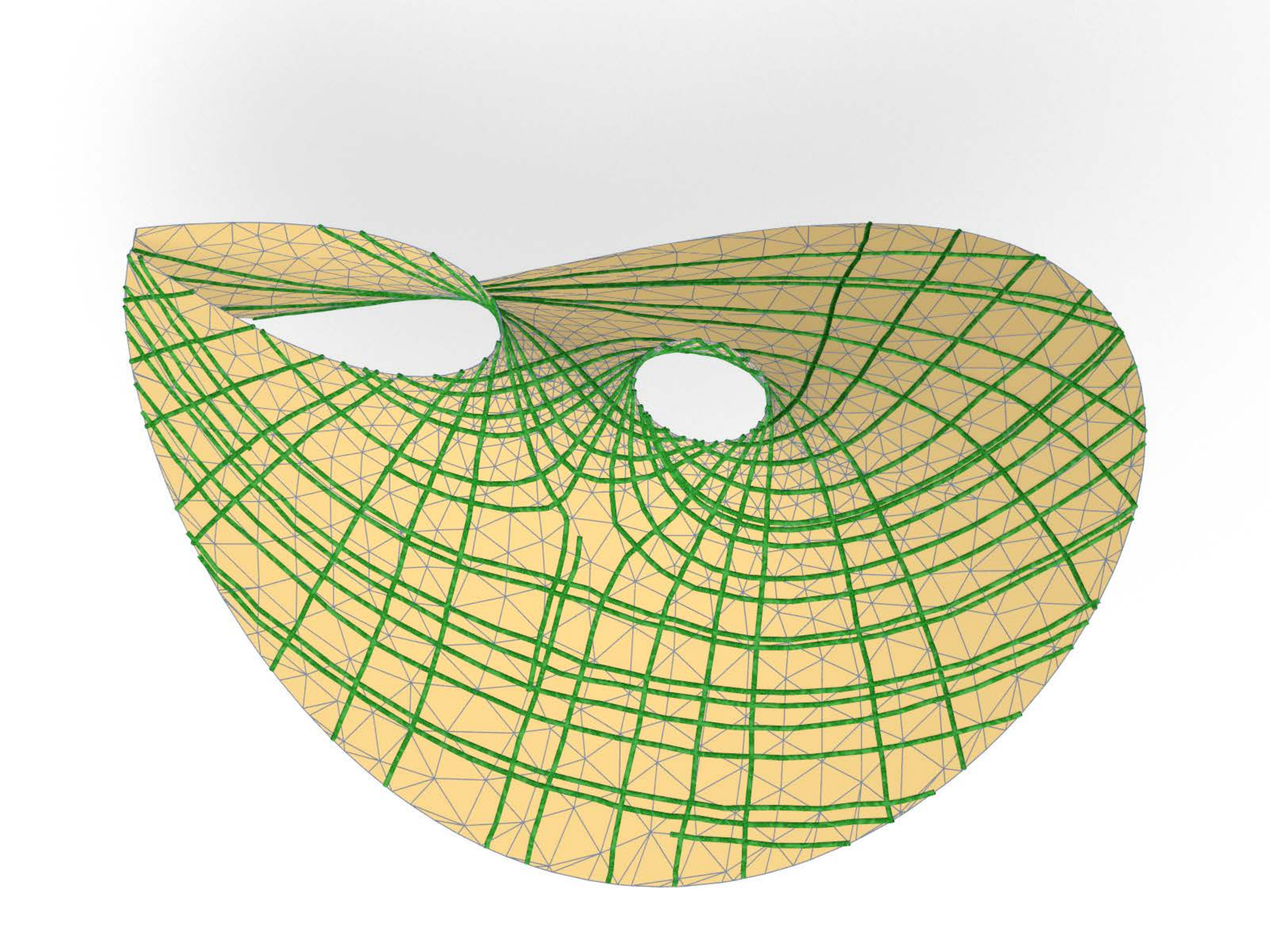}
     \subcaption{}
    \end{subfigure}
	\begin{subfigure}[t]{0.31\textwidth}
    \includegraphics[width=\textwidth]{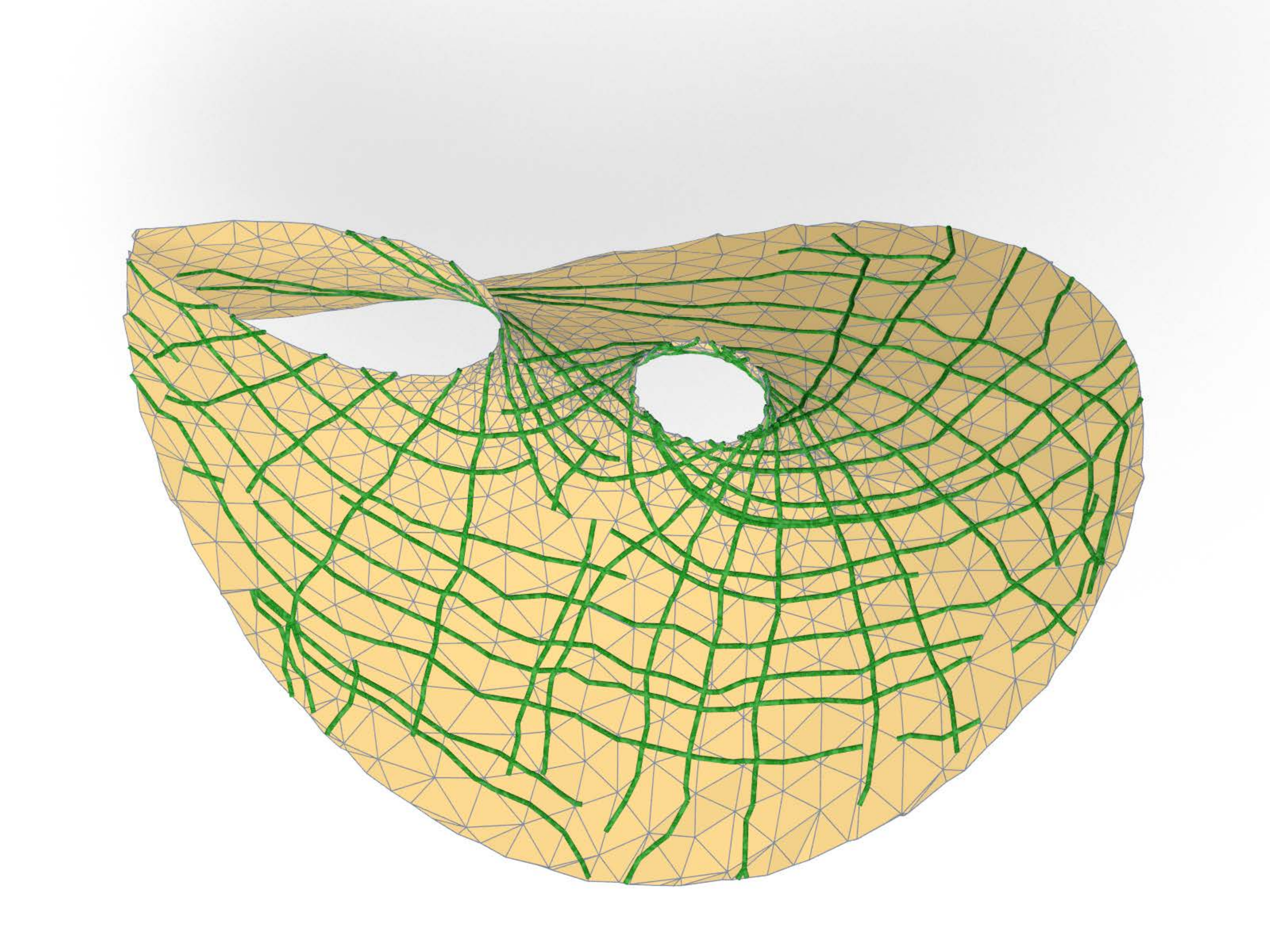}
    \subcaption{}\label{fig::Chennoisyr3}
    \end{subfigure}
    \begin{subfigure}[t]{0.31\textwidth}
    \includegraphics[width=\textwidth]{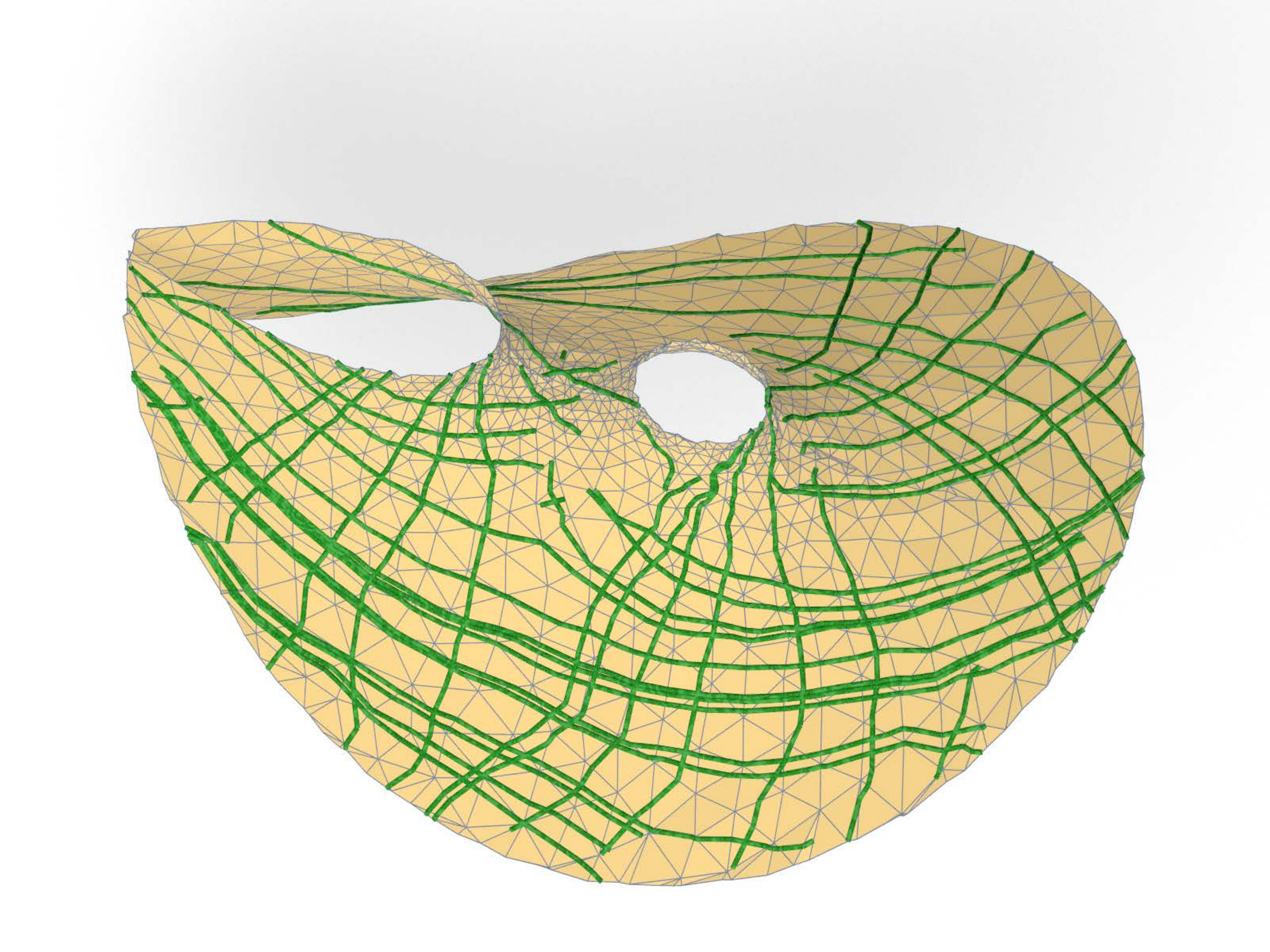}
    \subcaption{} \label{fig::Chennoisyr6}
    \end{subfigure}
\caption{By perturbing the positions of the vertices of the approximated Chen's surface, we get noisy data and corresponding noisy discrete asymptotic lines (a) The initial triangulation with asymptotic lines ($R=3$) (b) Asymptotic lines with noisy data and $R=3$ (c) Asymptotic lines with noisy data and $R=6$}
\end{figure}

\section{Conclusion and perspective}
 In this article, with any smooth subspace or singular geometric subspace $W$ of a Riemannian manifold $M$, we associate a family of cones,  defined over  any Borel subset of  $M$. These cones are the generalization of the asymptotic directions defined at each point of a smooth surface of the Euclidean space. We obtain convergence and approximation theorems, when a sequence of polyhedra tends to a smooth subspace. As a consequence, we find good approximations of the asymptotic lines of a surface when it is approximated by a suitable triangulation. In our future work, we will study the relations between the characteristic of these cones and the geometry (and topology) of $W$. We will apply these results in different fields (face recognition to detect similarities between scanned faces for instance).

\section*{Acknowledgments}
 We thank Fran\c{c}ois Golse and Simon Masnou for highlighting  interesting results  in measure theory that have been useful in our context, and Helmut Pottmann for his help and  judicious remarks on a first version of the text.

\bibliographystyle{plain}
\bibliography{mybibfile}

\end{document}